\documentclass[12pt]{article}
\usepackage[utf8]{inputenc}
\usepackage{amsmath, amsfonts, amscd, latexsym, amsthm, amssymb, graphicx}
\newtheorem{theor}{\hspace{1cm}{\sc Theorem}}[section]
\newtheorem{utver}[theor]{\hspace{1cm}{\sc Proposition}}
\newtheorem{sledst}[theor]{\hspace{1cm}{\sc Corollary}}
\newtheorem{lemma}[theor]{\hspace{1cm}{\sc Lemma}}
\newtheorem{conj}[theor]{\hspace{1cm}{\sc Conjecture}}
\newtheorem*{utver*}{\hspace{1cm}{\sc Proposition}}
\theoremstyle{definition}
\newtheorem{defin}[theor]{\hspace{1cm}{\sc Definition}}
\newtheorem{exa}[theor]{\hspace{1cm}{\sc Example}}
\newtheorem*{rem}{\hspace{1cm}{\sc Remark}}

\newcommand{\Vol}{\mathop{\rm Vol}\nolimits}

\newcommand{\codim}{\mathop{\rm codim}\nolimits}

\newcommand{\im}{\mathop{\rm im}\nolimits}

\newcommand{\sing}{\mathop{\rm sing}\nolimits}

\newcommand{\conv}{\mathop{\rm conv}\nolimits}

\newcommand{\dist}{\mathop{\rm dist}\nolimits}

\newcommand{\MV}{\mathop{\rm MV}\nolimits}

\newcommand{\MP}{\mathop{\rm MP}\nolimits}

\def\R{\mathbb R}
\def\N{\mathbb N}
\def\Z{\mathbb Z}
\def\Q{\mathbb Q}
\def\C{\mathbb C}
\def\CC{({\mathbb C}\setminus 0)}

\def\CP{\mathbb C\mathbb P}

\begin{document}
\title{Discriminant of system of equations}
\author{A. Esterov}
\date{}
\maketitle
{\begin{quote}
What polynomial in the coefficients of a system of algebraic equations should be called its discriminant?
We prove a package of facts that provide a possible answer. Let us call a system typical, if the homeomorphic type
of its set of solutions does not change as we perturb its (non-zero) coefficients. The set of all atypical
systems turns out to be a hypersurface in the space of all systems of $k$ equations in $n\ge k-1$ variables, whose
monomials are contained in $k$ given finite sets. This hypersurface $B$ is the union
of two well-known strata: the set of all systems that have a singular solution (this stratum is conventionally
called the discriminant) and the set of all systems, whose principal part is degenerate (they can be regarded
as systems with a singular solution at infinity). None of these two strata is a hypersurface in general, and
codimensions of their components have not been fully understood yet (e.g. dual defect toric varieties are not classified),
so the purity of dimension of their union seems somewhat surprising. We deduce it from a similar tropical purity fact
of independent interest: the stable intersection of a tropical fan with a boundary of a polytope in the ambient
space has pure codimension one in this tropical fan.

A generic system of equations in an irreducible component $B_i$ of the hypersurface $B$ always differs from a typical system
by the Euler characteristic of its set of solutions. Regarding the difference of these two Euler characteristics
as the multiplicity of $B_i$, we turn $B$ into an effective divisor, whose equation we call the
Euler discriminant of a system of equations by the following reasons. Firstly, it vanishes exactly at those systems
that have a singular solution (possibly at infinity). Secondly, despite its topological definition, it admits a
simple linear-algebraic formula for its computation, and a positive formula
for its Newton polytope. 
Thirdly, it interpolates many classical objects and inherits many of their nice properties: for $k=n+1$, it is
the sparse resultant (defined by vanishing on consistent systems of equations); for $k=1$, it is the principal
$A$-determinant (defined as the sparse resultant of the polynomial and its partial derivatives); as we specialize
the indeterminate coefficients of our system to be polynomials of a new parameter, the Euler discriminant turns
out to be preserved under this base change, similarly to discriminants of deformations. This allows, for example, to specialize
our results to generic polynomial maps: the bifurcation set of a dominant polynomial map, whose components are generic linear combinations of given monomials, is always
a hypersurface, and a generic atypical fiber of such a map differs 
from a typical one by its Euler characteristic.
\end{quote}}
{\scriptsize \tableofcontents}
\section{Introduction}
\subsection{Degenerate systems} For a finite set $H\subset\Z^n$, we study the space $\C[H]$ of all \textit{Laurent polynomials} $h(x)=\sum_{a\in H}c_{a}x^a$, where $x^a$ stands for the monomial $x_1^{a_1}\ldots x_n^{a_n}$, the coefficient $c_{a}$ is a complex number, and the polynomial $h$ is considered as a function $\CC^n\to\C$. For a linear function $v:\Z^n\to\Z$, denote the intersection of $H$ with the boundary of the affine half-space $H+\{v<0\}$ by $H^v$, and the highest $v$-degree component $\sum_{a\in H^v}c_ax^a$ by $h^v$ (if $v=0$, then we set $H^0=H$ and $h^0=h$).

In what follows, we denote a collection of finite sets $A_0,\ldots,A_k$ in $\Z^n$ by ${A}$, the space $\C[A_0]\oplus\ldots\oplus\C[A_k]$ by $\C[{A}]$, consider its element $f=(f_0,\ldots,f_k)\in\C[{A}]$ as a map $\CC^n\to\C^{k+1}$, and denote $(f_0^v,\ldots,f_k^v)$ by $f^v$.
\begin{theor} \label{thnondeg} 
Assume that $A_0+\ldots+A_k$ is not contained in an affine hyperplane. The following three conditions are equivalent for the system of equations $f=0$: \newline 1)
There exists an arbitrarily small $\tilde f\in\C[A]$, such that the sets $\{f=0\}$ and $\{f+\tilde f=0\}$ are not diffeomorphic.
\newline 2) There exists an arbitrarily small $\tilde f\in\C[A]$, such that the sets $\{f=0\}$ and $\{f+\tilde f=0\}$ have different Euler characteristic. \newline 3) There exists a linear function $v:\Z^n\to\Z$ such that the differentials $df_0^v,\ldots,df_k^v$ are linearly dependent at some point of the set $\{f^v=0\}$.
\end{theor}
This and subsequent theorems of the introduction are proved in Section \ref{Sunivcase}.
The assumption on $A_0+\ldots+A_k$ cannot be dropped, because, otherwise, the Euler characteristic of $\{f=0\}$ equals 0 for every $f\in\C[A]$ by homogeneity considerations.
\begin{defin} A system $f\in\C[A]$ is said to be \textit{degenerate}, if it satisfies any of the three conditions above.
\end{defin}
Condition 3 was introduced in \cite{kouchn} for $k=0$ and in \cite{kh77} for arbitrary $k$ in a slightly different from: for example, if $k=0$, then the (convex hull of) $A$ is assumed to be the Newton polytope of $f$ in the setting of \cite{kouchn}, while, in our setting, the Newton polytope of $f$ is contained in $A$, and $f$ is degenerate if its Newton polytope is strictly smaller than $A$. 

Condition 2 will play the role of tameness on a complex torus for our purpose (cf. the definition in \cite{brought}), although it is not equivalent to tameness at all. Its similarity
to tameness admits further development: e. g. for non-degenerate $f$ and a generic local system $L$ on $\CC^n$, so that $H(\CC^n,L)=0$, the twisted homology $H(\{f=0\},L)$ vanish except for the middle dimension. 

For example, if $k=n$, then $B$ is the resultant set (i.e. the set of all consistent systems of equations in $\C[A]$, see \cite{sturmf}); if $A_0=\ldots=A_k$ is the set of vertices of the standard $n$-dimensional simplex, then $f_0,\ldots,f_k$ are linear, and $B$ is defined by vanishing of the product of the maximal minors for the matrix of coefficients of $f_0,\ldots,f_k$.
\begin{defin} \label{defrel} The collection $A$ is said to be \textit{relevant}, if the dimension of the convex hull of $A_{i_0}+\ldots+A_{i_p}$ is at least $p$ for every sequence $0\leq i_0<\ldots<i_p\leq k$, and equals $n$ for $p=k$.
\end{defin}
\begin{theor} \label{lozung0} If $A$ is relevant, then the set $B$ of all degenerate systems in $\C[A]$ is a non-empty hypersurface.
\end{theor}
The assumption of relevance cannot be dropped, because, otherwise, the set of consistent systems has codimension greater than 1 (see \cite{sturmf}). 

The similar question of whether the $A$-discriminant $\{ f\in\C[A]\, |\, f=0 $ is not regular$\}$ is a hypersurface is well known for $k=0$ as the problem of classification of dual defect polytopes, and is still open (see \cite{cat}, \cite{roc}, \cite{dfs}, \cite{tak}, \cite{dcg}, etc). Moreover, for $k>0$, the $A$-discriminant may be not of pure dimension: e.g. for $A_0=\{0,1\}\times\{0,1\}$ and $A_1=\{0,1,2\}\times\{0\}$, there is a codimension 1 component, to which $f=\Bigl(a+bx+cy+dxy,\, r(x-p)^2\Bigr)$ belongs, and a codimension 2 component, to which $f=\Bigl(a(x-b)(y-c),\, r(x-b)(x-p)\Bigr)$ belongs. ``Fortunately'', the latter one is swallowed up by the codimension 1 stratum of $B$, to which $f=\Bigl(a(x-b)(y-c)+d,\, r(x-b)(x-p)\Bigr)$ belongs because of its singularity at infinity. The generic configuration of $f_0=0$ (in solid lines) and $f_1=0$ (in dotted lines) is shown on the picture below, followed by the configurations of the three mentioned degenerations.

\noindent\includegraphics[width=\textwidth]{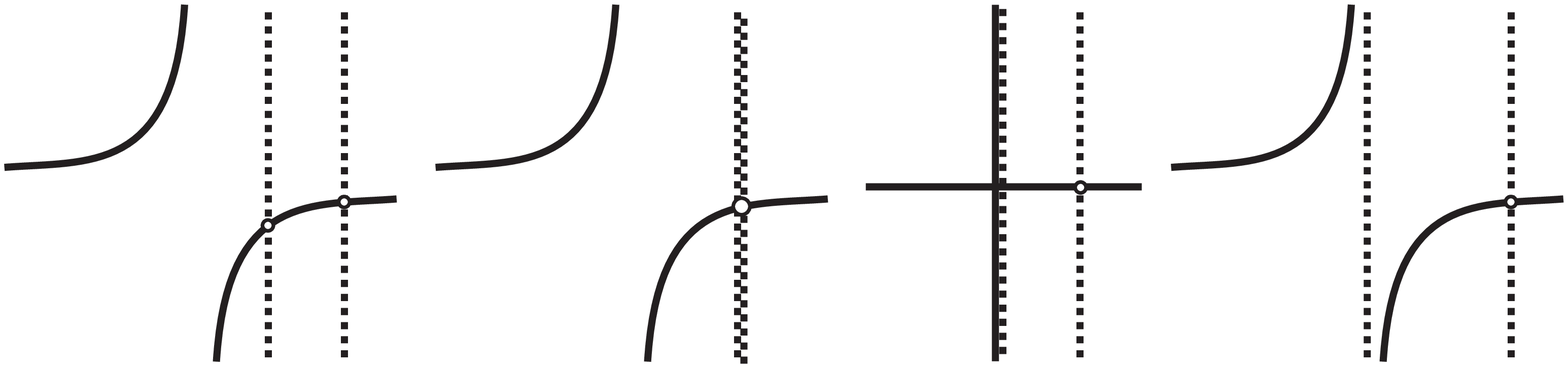}

If, instead of degenerate Laurent polynomials on $\CC^n$ we study degenerate polynomials on $\C^n$ (i. e. polynomials $f\in\C[A]$ such that $\{x\in\C^n\, |\, f(x)=0\}$ is not homeomorphic to $\{x\in\C^n\, |\, \tilde f(x)=0\}$ for a close $\tilde f\in\C[A]$), then the set $B$ may be of higher codimension, as it happens for $n=2,\, k=0, \, A=\{(0,0), (0,1), (1,0)\}$, and it would be interesting to classify all such exceptional collections $A_0,\ldots,A_k$ in $\Z^n_{\geqslant 0}$, containing $0$.

\subsection{Euler discriminant}
For $f\in\C[A]\setminus B$ and generic $\tilde f$ in an irreducible component $B_i\subset B$, denote the difference of the Euler characteristics $e\{\tilde f=0\}-e\{f=0\}$ by $e_i$.
\begin{utver} \label{coinc0} If $A$ is relevant, then $e_i>0$ for $n-k$ even and $e_i<0$ for $n-k$ odd.
\end{utver} For the proof, see Corollary \ref{coinc1}.
\begin{defin} \label{defaeuldiscr} If $A$ is relevant, then the equation of the effective divisor $(-1)^{n-k}\sum_ie_iB_i$ is called the $A$-\textit{Euler discriminant} and is denoted by $E_A=E_{A_0,\ldots,A_k}$.
\end{defin}
By Proposition \ref{coinc0}, $E_A$ is a non-constant polynomial on $\C[A]$, defined up to multiplication by a non-zero constant. By Theorem \ref{thnondeg}, the equation $E_A=0$ describes all degenerate systems of equations in $\C[A]$. It can be computed as follows.

Recall that, for a finite subset $H$ in an affine hyperplane of $\Z^m$, the $H$-\textit{sparse resultant} $R^{red}_H$ is the unique (up to multiplication by a constant) polynomial on $\oplus^{m}\C[H]$, vanishing at all points $(h_1,\ldots,h_m)\in\oplus^{m}\C[H]$, such that the (quasihomogeneous) system $h_1=\ldots=h_m=0$ has a solution in $\CC^m$ (see \cite{sturmf} for details).  The $H$-principal determinant $E^{red}_H$ is defined at $f\in\C[H]$ as $R^{red}_H(x_1\frac{\partial f}{\partial x_1},\ldots,x_m\frac{\partial f}{\partial x_m})$, where $x_1,\ldots,x_m$ are the standard coordinates on $\CC^m$ (see \cite{gkz} for details).

Let $\lambda_0,\ldots,\lambda_k$ be the standard coordinates on $\C^{k+1}$, and let $e_0,\ldots,e_k$ be the standard basis in the lattice $\Z^{k+1}$. For every $I\subset\{0,\ldots,k\}$, denote the lattice with the generators $\{e_i\, |\, i\in I\}$ by $\Z^I$, the \textit{Cayley configuration} $\bigcup_{i\in I} \{e_{i}\}\times A_{i}\subset\Z^I\times\Z^n$ by $A_I$, and the polynomial $\sum_{i\in I}\lambda_if_i\in\C[A_I]$ 
by $f_I$. 

\begin{theor} 
\label{lozungdecomp} 
We have $$E_A(f)=\prod_{I\subset\{0,\ldots,k\}}E^{red}_{A_I}(f_I)^{(-1)^{k+1-|I|}\left|\frac{\Z^I\times\Z^n}{A_I}\right|},$$ 
where the product is taken over all $I$ such that $A_I$ generates an $(n+|I|)$-dimensional lattice.
\end{theor} 
This reduces the computation of the Euler discriminant to the computation of sparse resultants or principal determinants, which are well studied (see e.g. \cite{gkz} and \cite{sturmf}).
In particular, if $k=n$, then $E_A$ is the sparse resultant, and for $k=0$ it coincides with the $A$-principal determinant. Similarly to these special cases, the Newton polytope of the Euler discriminant admits the following simple description.

\subsection{Degree of Euler discriminant}
Recall that a coherent triangulation $T$ of a finite subset $H\subset\Z^N$ is a set of $N$-dimensional simplices, such that

1) their vertices are contained in $H$,

2) the union of them is the convex hull of $H$,

3) the intersection of any two of them is their common face (maybe empty),

4) they are the domains of linearity of a convex piecewise-linear function (this property is called {\it coherence} or {\it convexity}).

Coefficients of polynomials in $\C[A]$ form a natural coordinate system $(c_{a,i})_{a\in A_i,\, i=0\ldots,k}$ on it, so that $c_{a,i}(f)$ is the coefficient of the monomial $x^a$ in the $i$-th polynomial of the tuple $f\in\C[A]$.
For a simplex $S$ with vertices in  $A_{\{0,\ldots,k\}}$, let $c_S$ be the product of all $c_{a,i}$, such that $(e_i,a)$ is a vertex of $S$, and, for every $j\ne i$, 
the set $\{e_j\}\times \Z^n$ contains more than one vertex of $S$.

\begin{utver} \label{lozungnewt}
The set of monomials 
$$\bigl\{\prod_{S\in T}c_S^{\Vol S}\, |\, T \mbox{ is a coherent triangulation of } A_{\{0,\ldots,k\}}\bigr\}$$ 
is the set of vertices for the Newton polytope of the $A$-Euler discriminant $E_A$ (the Newton polytope is in the natural coordinate system $(c_{a,i})$, and $\Vol$ stands for the integer volume, normalized by the condition $\Vol($standard simplex$)=1$).
\end{utver}
This follows from the description of the Newton polytope of the Euler discriminant (Theorem \ref{thmain}.3) and of the vertices of a mixed fiber polytope (Theorem 3.19 in \cite{ekh}). 

In particular, the degree of the Euler discriminant is as follows (the same formula can be deduced from Theorem \ref{lozungdecomp}, which reduces it to the known special case of $A$-principal determinants). 

Recall that the {\it mixed volume} of finite sets $B_1,\ldots,B_n$ in $\R^n$ is the symmetric function, taking values at $\R$, additive with respect to the Minkowski summation $P+Q=\{p+q\, |\, p\in P,\, q\in Q\}$,
and equal to the volume of the convex hull $\conv B$ for $B_1=\ldots=B_n=B$. We denote it by the monomial $B_1\cdot\ldots\cdot B_n$. In the same way, the {\it mixed moment} of finite sets $B_0,\ldots,B_n$ in $\R^n$ is the symmetric function, taking values at $\R^n$, additive with respect to the Minkowski summation,
and equal to the moment $\int_{\conv B}x\, {\rm d}x$ for $B_0=\ldots=B_n=B$. We denote it by the monomial $B_0\cdot\ldots\cdot B_n$. As usual, the formal derivative $\frac{\partial}{\partial B_0}(B_0^{a_0}\cdot\ldots\cdot B_n^{a_n})$  is defined as $a_0 B_0^{a_0-1}\cdot B_1^{a_1}\cdot\ldots\cdot B_n^{a_n}$.

\begin{sledst} \label{lozungdeg} 1) As a polynomial of the variables $c_{a,i},\, a\in A_i$, for a given $i$, the $A$-Euler discriminant is homogeneous of degree $$\deg_{n,i} A=\frac{\partial}{\partial A_i} \left(\sum_{(a_0,\ldots,a_k)\in\N^{k+1}, \atop a_0+\ldots+a_k=n+1}  A_0^{a_0}\cdot\ldots\cdot A_k^{a_k}\right).$$

2) For a linear function $v:\Z^n\to\Z$, define the $v$-quasidegree of the variable $c_{a,i}$ as $v(a)$. Then the $A$-Euler discriminant is $v$-quasihomogeneous of degree $$\deg_{n,v} A=v \left(\sum_{(a_0,\ldots,a_k)\in\N^{k+1}, \atop a_0+\ldots+a_k=n+1}  A_0^{a_0}\cdot\ldots\cdot A_k^{a_k}\right).$$

3) Every quasihomogeneity of the $A$-Euler discriminant is a linear combination of the homogeneities, described in (1) and (2).
\end{sledst}

Proposition \ref{lozungnewt} implies that the Newton polytope of the Euler discriminant is the mixed secondary polytope of $A$ in the sense of the following definition (well known for $k=0$).
\begin{defin} For a coherent triangulation $T$ of finite sets $A_0,\ldots,A_k$ in $\R^k$, let $T_{i,a}$ be the sum of the integer volumes of all simplices $S\in T$,
such that $(e_i,a)$ is a vertex of $S$ and, for every $j\ne i$, the set $\{e_j\}\times \Z^n$ contains more than one vertex of $S$. The convex hull of the set of points with coordinates $(T_{i,a}),\, i=0,\ldots,k,\, a\in A_i,$ in the lattice $\Z^{A_0}\oplus\ldots\oplus\Z^{A_k}$, as $T$ runs over all coherent triangulations of $A_0,\ldots,A_k$, is called the {\it mixed secondary polytope} of $A_0,\ldots,A_k$ and is denoted by $S_A$.
\end{defin}
Note that,
denoting the vertices of $S$ in $\{e_j\}\times \R^n$ by $s^i_0,\ldots,s^i_{n_i}$, we have $n_0+\ldots+n_k=n$, and the volume of $S$ equals the determinant of the $n\times n$ matrix, whose rows are $s^i_j-s^i_0$ for $i=0,\ldots,k$ and $j=1,\ldots,n_i$.

\subsection{Factorization of Euler discriminant}
Let us describe the prime factorization of the Euler discriminant $E_A$. We denote the Gelfand-Kapranov-Zelevinsky $A_{\{0,\ldots,k\}}$-discriminant (which is the equation of the set $\{F\, |\, F=0$ is not regular$\}\subset\C[A_{\{0,\ldots,k\}}]$, see \cite{gkz}) by $\Delta_A$.
(This notation agrees with the one for mixed discriminants in \cite{st11}, see Theorem 2.1 therein.) The polynomial $\Delta_{A}$ is prime, and Theorem \ref{lozung0} and Corollary \ref{coinc1} give the following decomposition of the Euler discriminant into such polynomials. A collection $\Gamma$ of non-empty subsets $\Gamma_i\subset A_i,\, i\in I\subset\{0,\ldots,k\}$, is called a {\it facing} of $A$, if $\Gamma_I\subset\Z^I\times\Z^n$ equals $A^v_{\{0,\ldots,k\}}$ for some linear function $v:\Z^{k+1}\times\Z^n\to\Z$. In this case, we denote the polynomial $f^v_{\{0,\ldots,k\}}$ by $f_\Gamma$, and the index of the lattice, generated by the pairwise differences of the points of $\Gamma_I$, by $i_\Gamma$. 
\begin{utver} We have $$E_A(f)=\prod_{\Gamma}\Delta_{\Gamma}(f_{\Gamma})^{i_\Gamma c_A^\Gamma},$$ where $\Gamma$ runs over all facings of $A$, and $c_A^\Gamma$ is the Milnor number of $A$ at $\Gamma$ (see Definition \ref{mixmilnum}).
\end{utver}
Passing to the Newton polytopes of both sides of this equality, we obtain the following representation for the mixed secondary polytope:
$$S_A=\sum_{\Gamma} i_\Gamma c_A^\Gamma \cdot(\mbox{Newton polytope of } \Delta_A),\eqno{(*)}$$
where the summation is the Minkowski summation $P+Q=\{p+q\, |\, p\in P,\, q\in Q\}$ and the multiplication is the dilatation. Since the mixed secondary polytope is defined combinatorially, while the Newton polytope of $\Delta_A$ is not, the inverse of this formula is more useful.
For $\Gamma$ and $B$ running over all facings of $A$, let $(e^\Gamma_B)$ be the inverse to the matrix $(i_\Gamma c^\Gamma_B)$, where by convention $c^\Gamma_B=0$ unless $\Gamma$ is a facing of $B$. Note that the matrix $(i_\Gamma c^\Gamma_B)$ is upper triangular, nondegenerate and quite sparse (see Corollary \ref{importpositive}).
\begin{defin} The number $e^\Gamma_A$ is called the {\it Euler obstruction} of $A$ at $\Gamma$.
\end{defin}
If $k=0$, then it is indeed the Euler obstruction of the $A$-toric variety at a point of its $\Gamma$-orbit, see \cite{tak} or \cite{dcg} for a precise statement. Inverting the system of linear equations $(*)$ for all $A$ and $\Gamma$, we obtain the following.
\begin{sledst} The Newton polytope of $\Delta_A$ equals
$$\sum_{\Gamma} e^\Gamma_A \cdot S_\Gamma.$$
\end{sledst}
This formula specialzes to the Gelfand-Kapranov-Zelevinsky formula for the Newton polytope of the $A$-discriminant (if $k=0$, \cite{gkz}, \cite{tak}), to similar formulas
for the Newton polytopes of the $A$-resultant (if $k=n$, \cite{sturmf}) and the mixed discriminant (if $k=n-1$, \cite{st11}), and for the degree of the classical discriminant in terms of the degree of $f_i$ (if $A_i$ are multiples of the standard simplex, see \cite{benoist} or \cite{nie}). Note that the formula is not positive, because Euler obstructions may be negative; see \cite{dcg} for a discussion of Minkowski subtraction of polytopes in this context, and see \cite{dfs} and \cite{st11} for a positive formula, based on a very different idea.

In particular, the degree of $\Delta_A$ can be computed by Corollary \ref{lozungdeg} as follows. For a facing $\Gamma$, we denote the dimension of the convex hull of $\sum\Gamma_i$ by $\dim_\Gamma$.
\begin{sledst} \label{lozungdegirred} 1) As a polynomial of the coefficients of $f_i$, the discriminant $\Delta_A(f_0,\ldots,f_k)$ is homogeneous of degree $$\sum_{\Gamma} e^\Gamma_A \cdot \deg_{\dim_\Gamma,i}\Gamma,$$
where the summation is taken over all facings $\Gamma$ of $A$, whose elements are indexed by a set, containing $i$.

2) For a linear function $v:\Z^n\to\Z$, the polynomial $\Delta_A$ is $v$-quasihomogeneous of degree $$\sum_{\Gamma} e^\Gamma_A \cdot \deg_{\dim_\Gamma,v}\Gamma,$$
where the summation is taken over all facings $\Gamma$ of $A$.

See Corollary \ref{lozungdeg} for the meaning of $\deg_{\cdot,i}$, $\deg_{\cdot,v}$ and $v$-quasihomogeneity.\end{sledst}

Note that, for some $A$, there exist quasihomogeneities of $\Delta_A$ that are not linear combinations of the ones listed above. It would be a step towards the classification of dual defect collections $A$ to classify those with ``many'' additional quasihomogeneities of $\Delta_A$ (or $D^{red}_A$, introduced in \cite{dcg}).

\subsection{Relation to bifurcation sets} \label{Sintrobif} We prove the formulated results in a more general setting that also covers a number of other noteworthy special cases. Namely, we assume the coefficients $c_{a,i}$ of each polynomial $f_i(x)=\sum_{a\in A_i}c_{a,i}x^a$ to be not independent variables, but polynomials on a parameter space $M$. We consider two different cases.
\newline \textsc{I (Laurent polynomials).} We have $M=\CC^m$, pick a finite set $H_{a,i}\subset\Z^m$ for every $a\in A_i$, and consider $c_{a,i}\in\C[H_{a,i}]$ as a function on $M$.
\newline \textsc{II (polynomials).} We have $M=\C^m$, pick a finite set $H_{a,i}\subset\Z^m_{\geqslant 0},\, 0\in H_{a,i},$ for every $a\in A_i$, and consider $c_{a,i}\in\C[H_{a,i}]$ as a function on $M$.

In both cases, we denote the union $\bigcup_{a\in A_i} \{a\}\times H_{a,i}\subset\Z^n\oplus\Z^m$ by $H_i$, so that $f_i\in\C[H_i]$ is a function on $\CC^n\times M$. As before, we abbreviate the collection $H_0,\ldots,H_k$ to $H$, and the vector polynomial $(f_0,\ldots,f_k)$ to $f$. Seeing $f(x,t)$ as a vector function of $x\in\CC^n$ with fixed $t\in M$, we consider it as an element of $\C[A]$ and denote it by $F(t)$, so that $F$ is a polynomial map $M\to\C[A]$.

We study the projection $\pi$ of the set $Q=\{f=0\}\subset\CC^n\times M$ to $M$.
\begin{defin} \label{defbif0} The \textit{bifurcation set} $B_p$ of a morphism of algebraic varieties $p:Q\to M$ is the complement to the maximal open set $S\subset M$, such that the restriction of $p$ to the preimage $p^{-1}(S)$ is a locally trivial fibration.
\end{defin}
\begin{theor} \label{lozung} Assume that $A$ is relevant, $H$, $M$ and $F$ are defined as in one of the settings {\rm(I)} or {\rm(II)} above, and $f\in\C[H]$ is nondegenerate (or, less restrictively, $f$ is 1-nondegenerate in the sense of Definition \ref{cdeg} below). \newline
1) The bifurcation set $B_\pi$ is given by the equation $\{E_A\circ F=0\}$. In particular, $B_\pi$ is a hypersurface (maybe empty). \newline 2) For generic points $t_0\in B_\pi$ and $t\notin B_\pi$, the local multiplicity of $E_A\circ F$ at $t_0$ equals $(-1)^{n-k}$ times the difference of the Euler characteristics of the fibers $$e\{f(\cdot,t_0)=0\}-e\{f(\cdot,t)=0\}.$$ 
\newline 3) The bifurcation set is empty if and only if
$$\dim(H_{i_1}+\ldots+H_{i_p})<p$$ for some $i_1<\ldots<i_p$, or
$$\dim(H_0+\ldots+H_k)\leqslant n.$$
\end{theor} For the proof of Parts 1 and 2, see Theorem \ref{thmain}.1 and 2. Part 3 follows from Proposition \ref{nonempty1} and describes the two obvious cases of an empty bifurcation set: in the first case, $Q$ is empty, and, in the second case, $Q=Q'\times M$ for some $Q'\subset\CC^n$ up to a monomial change of coordinates. Part 2 implies that the equation $E_A\circ F=0$ for the set $B_\pi$ is not square-free, but the multiplicities of its factors measure degeneration of fibers of $\pi$. Part 1 is a manifestation of a more general phenomenon of ``invariance of the bifurcation set under base changes'', see e.g. \cite{bodin} for a similar phenomenon in two variables.
Theorem \ref{lozung} proves the following conjecture in the special case of nondegenerate complete intersections in the torus.
\begin{conj} For the projection $p$ of an arbitrary sch\"{o}n submanifold of a complex torus (see \cite{tev}) along arbitrary subtorus, the bifurcation set $B_p$ has codimension 1, and a generic fiber of $p$ over $B_p$ differs from a fiber of $p$ outside $B_p$ by its vector of virtual Betti numbers (and, a forteriori, by its Hodge-Deligne polynomial).
\end{conj}
We can relax the assumption of nondegeneracy in Theorem \ref{lozung} 
as follows. The \textit{dimension} of a linear function $v:\Z^n\oplus\Z^m\to\Z$ is the maximal dimension of a fiber for the projection of the convex hull of $\sum_i H_i^v\subset\Z^n\oplus\Z^m$ to $\Z^n$.
\begin{defin} \label{cdeg} The vector polynomial $f\in\C[H]$ is said to be $i$-\textit{degenerate}, if there exists a linear function $v:\Z^n\oplus\Z^m\to\Z$ of dimension at most $i$, such that the differentials $df_0^v,\ldots,df_k^v$ are linearly dependent at some point of the set $\{f^v=0\}$. 
\end{defin}
The nondegeneracy assumption of Theorem \ref{lozung} can be relaxed 
to 1-nondegeneracy, but not to 0-nondegeneracy: for example, the bifurcation set of the projection of $1+x+z+yz=1+x+bz+byz=0$ to the $(x,y)$-plane is $\{(0,0)\}$ for $b\ne 1$, in spite of 0-nondegeneracy of these equations. There is also a version of Theorem \ref{lozung} for irrelevant $A$.
\begin{utver} \label{lozungirr} Let $c$ be the maximum of the differences $$p-\dim(\mbox{convex hull of } A_{i_1}+\ldots+A_{i_p})$$ over all sequences $0\leq i_1<\ldots<i_p\leq k$. If $c>0$, and $f$ is $c$-nondegenerate, then the bifurcation set $B_\pi$ equals the closure of the image of $\pi$ and has pure codimension $c$. Moreover, a generic atypical fiber of $\pi$ has a nonzero Euler characteristic, provided that $A_0+\ldots+A_k$ is not contained in an affine hyperplane.
\end{utver} For the proof, see Proposition \ref{irrelpure}.
Note that, for $\min\dim A<-1$, one can also study the bifurcation set of the projection $\pi:Q\to \pi(Q)$; this is a different problem, which makes sense and will be addressed in a subsequent paper.

The number $c$ in the statement of Proposition \ref{lozungirr} is called the {\it dimension of} the collection $A$.

\subsection{Application to topology of polynomial maps}
For a nondegenerate tuple of polynomials $g\in\C[A],\, g:\CC^n\to\C^{k+1}$,
we can apply Theorem \ref{lozung} in the setting (II) to the graph $g_0(x)-\lambda_0-c_0=\ldots=g_k-\lambda_k-c_k=0$ in $\CC^n\times\C^{k+1}$, shifted by generic numbers $c_i\in\C$, and get the following.
\begin{sledst} \label{lozung2} If $A$ is relevant, the tuple $g\in\C[A]$ is nondegenerate, and every $k$ of the $k+1$ polynomials in $g$ also form a nondegenerate tuple, then the bifurcation set of the map $g:\CC^n\to\C^{k+1}$ is a hypersurface, and a generic atypical fiber of $g$ differs from a typical one by its Euler characteristic.
\end{sledst}
Theorem \ref{thmain}.3 describes the Newton polytope of this hypersurface. By Theorems \ref{thdcg}, \ref{lozung}.1, \ref{theuler} and Corollary \ref{importpositive}, the difference of the Euler characteristics of the fibers in this statement consists of contributions of singularities (in the domain or at infinity) of the atypical fiber of $g$: every its singularity corresponds to a certain \textit{important face} $\Gamma$ of $A$ (see Definition \ref{defimport}), and its contribution to the difference of the Euler characteristics equals the \textit{Milnor number} of $A$ at $\Gamma$ (see Definition \ref{mixmilnum}). We prefer not to make this statement precise here, but rather refer to Section \ref{latticesets} for the complete picture, because details are subtle (e. g. the aforementioned singularities at infinity may be non-isolated even in the minimal toric compactification of $g$, compatible with $A$).

One part of Corollary \ref{lozung2} addresses the question of purity of the bifurcation locus, which is trivial for $k=0$ and is classical for $k=n-1$: see \cite{j1} and \cite{j2} for the purity of the Jelonek set.

Another part of this statement addresses the question of distinguishing atypical fibers of polynomial maps by their discrete invariants. This question, in contrast to the first one, is trivial for $k=n-1$ and is classical for $k=0$: see e. g. \cite{vl} and \cite{v1} for the case of two variables, \cite{stduke}, \cite{par1}, \cite{par2}, \cite{alm} and \cite{alm2} over $\C$  
and \cite{raibaut} over an arbitrary field of characteristic 0 for polynomials with isolated singularities at infinity, \cite{nz} and \cite{z2} for non-degenerate polynomials. This question was also addressed for $k=n-2$ (see e. g. \cite{vt1} and \cite{vt2}), but less is known for arbitrary $k$ (see e. g. \cite{tib} and \cite{gaf}). Note that the most common setting for these studies is the assumption of isolated singularities at infinity, which is neither weaker nor stronger than the one in Corollary \ref{lozung2}. Polynomials with isolated singularities at infinity may be degenerate (e.g. $(x-y)^2+(x-y)+c$), and, vice versa, nondegenerate polynomials, whose Newton polytopes are not simple, may have non-isolated singularities in any smooth compactification (e.g. $xyz(x+y)(z+1)+c$). It is an interesting problem of toric singularity theory to unify these two settings.

Many of the aforementioned works are also concerned with estimating the degree of the bifurcation set in terms of the degree of the mapping (see also \cite{lo1} and \cite{lo2} for $k=0, n=2$, \cite{j31} and \cite{jk1} for $k=0$, and \cite{j32} and \cite{jk2} for the general case). For nondegenerate maps (or, more generally, for 1-nondegenerate maps, see Definition \ref{cdeg}), the precise answer regarding the degree is given by Theorem \ref{thmain}.3 (it, moreover, describes the Newton polytope). For 0-nondegenerate maps, the bifurcation set may have components of higher codimension, but still it is contained in a hypersurface, whose degree is given by the formula of Theorem \ref{thmain}.3. For general $n$ and $k$, one cannot relax the assumption of 0-nondegeneracy, although, for some extreme cases, parts of this statement remain valid with no assumption on genericity of the mapping (see \cite{j1} for the purity of the ramification locus of an arbitrary generically finite map, and \cite{vl} for distinguishing atypical values of an arbitrary bivariate polynomial by their Euler characteristic).

The collection of references given in this subsection is very narrow and not intended to be an overview of the extensive literature on singularities at infinity. 

\subsection{Structure of paper} Section 2 prepares several necessary facts from tropical and convex geometry. All of these facts reflect the well known relation between linear dependence of polytopes, vanishing of their mixed volume, and emptiness of the stable intersection of their tropical fans (see Lemma \ref{lequiv1}). We also prove a refined version of purity for the dimension of the stable intersection of tropical varieties: the stable intersection of a $k$-dimensional affine tropical variety with the boundary of a polytope in the ambient space is always purely $(k-1)$-dimensional (Proposition \ref{dtrop}).
It will be used to prove purity of dimension for bifurcation sets over $\C$ and can be seen itself as purity of dimension for a certain tropical bifurcation set (see Conjecture \ref{tropjel}).

In Section 3, we estimate the bifurcation set of an arbitrary morphism $p$ of algebraic varieties from above (by the \textit{Bertini discriminant}, Definition \ref{defbertdiscr}) and from below (by the \textit{Euler discriminant}, Definition \ref{defeuldiscr}). The Bertini discriminant is the set of critical values of a stratified fiberwise compactification of $p$; it is widely used to estimate the bifurcation set from above. The Euler discriminant is the codimension 1 part of the discriminant in the sense of Teissier \cite{teissier}; the codimension 1 part makes sense for morphisms of algebraic manifolds with arbitrarily complicated singularities of fibers.

The Euler discriminant turns out to be invariant under base changes, subject to certain transversality assumptions (Theorem \ref{thbasech}). The Bertini discriminant turns out to be a hypersurface (Theorem \ref{thgeneral}.4), whenever the strata of the fiberwise compactification are  affine complete intersections, whose images are adjacent to the ones of codimension 1.

In Section 4, we apply our estimates of the bifurcation set to the setting of Theorem \ref{lozung}. All the results, mentioned in the introduction, follow from the three observations regarding Bertini and Euler discriminants in this toric setting:
\newline 1) The map $F$ that appears in Theorem \ref{lozung} satisfies the transversality assumptions of Theorem \ref{thbasech}, assuring the invariance of the Euler discriminant under the base change $F$, if and only if $f$ is $1$-nondegenerate.
\newline 2) The Bertini discriminant is a hypersurface, because a smooth toric compactification of $\CC^n$, compatible with $A$, satisfies the assumption of Theorem \ref{thgeneral}.4. If we have $A_0=\ldots=A_k$, then the adjacency condition of Theorem \ref{thgeneral}.4 reflects the obvious fact that every face of the convex hull of $A_0$ is a facet of another face. However, if $A_i$ are not equal to each other, then we need a "mixed" version of this obvious fact. Amazingly, this mixed version is the principal combinatorial challenge in our work, and its solution occupies most of Section 2 (see Theorem \ref{thess_main} for the precise statement). We reduce the problem to the aforementioned tropical fact about purity of stable intersections. As a result, the mixed version of the poset of faces of a polytope is understood well enough for our purpose, but not completely: for instance, the version that we construct is not a poset. See Conjecture \ref{conjess_main} for a possible way to fix it.
\newline 3) Both the Euler and the Bertini discriminant can be decomposed into irreducible components; an earlier paper \cite{dcg} was devoted to the study of these irreducible components, and its results imply that the decompositions of  the Euler and the Bertini discriminant coincide in the toric setting. Thus, in this setting, both the Euler and the Bertini discriminant coincide with the bifurcation set. That is why we do not discuss other known notions of the discriminant set of a non-proper map (e.g. the set of critical values at infinity in the sense of Ehresmann, see \cite{kurdyka} or \cite{rabier}): all such discriminants are also obviously estimated by the Euler and the Bertiny discriminant from below and from above respectively, and thus all of them coincide under our assumptions. 

At the end of Section 4, we assemble these three observations into the proof of the package of facts, presented in the introduction. 


I am very grateful to Eduardo Cattani, Alicia Dickenstein, Sandra Di Rocco, Pedro Gonzalez Perez, Askold Khovanskii, Alejandro Melle Hernandez, Benjamin Nill, Michel Raibaut, Martin Sombra, Bernd Sturmfels and Kiyoshi Takeuchi for enlightening and fruitful discussions. The text was thoroughly checked and greatly improved by Gleb Gusev; he generously provided a new proof of Theorem \ref{trivmp} and many other ideas that significantly simplify the final version of the paper. 



\section{Preliminaries from polyhedral geometry}
In Section \ref{Smixvolstabins}, we recall the notion of stable intersection and the relation between linear dependence of polytopes, vanishing of their mixed volume, and emptiness of the stable intersection of their dual complexes (Lemma \ref{lequiv1}). We also find a local version of this relation (where one of the polytopes is a cone, see Lemma \ref{lequiv_rel}) and a stronger version of purity for stable intersections of tropical complexes (with the boundary of a polytope instead of one of the tropical complexes, see Proposition \ref{dtrop}). In the toric setting, these two combinatorial facts will lead to the equality between the Euler and the Bertini discriminant and to the purity of the Bertini discriminant respectively.

In Section \ref{Smixfib}, we recall the notion and properties of mixed fiber polytopes and establish a version of Lemma \ref{lequiv1} for them, giving a criterion for their vanishing, see Theorem \ref{trivmp}. In order to formulate it, we recall the notion of an essential tuple of polytopes (Definition \ref{defess}). To keep computations traceable, we introduce certain vector notation for tuples of polytopes, and use it in the rest of Section 2.

In Section \ref{Sessfaces}, we look for a mixed version of an obvious fact: every face of a polytope is a face of its facet. In order to do so, we construct a mixed version of the poset of faces of a polytope, corresponding to a tuple of polytopes. This mixed version is not a poset (its order relaton is not transitive) but behaves nicely enough for our purpose, see Theorem \ref{thess_main}.

In Section \ref{latticesets}, we apply the aforementioned combinatorial results to study the set of degenerate systems of equations. In particular, we prove that a generic system from its codimension 1 part differs from a nondegenerate system by the Euler characteristic of the set of solutions, see Corollary \ref{coinc1}. We also suggest a possibly better candidate to the mixed version of the poset of faces, see Conjecture \ref{conjess_main}.

\subsection{Mixed volumes and stable intersections} \label{Smixvolstabins} A \textit{polyhedral complex}  of dimension $k$ in $\R^n$ is a locally finite union of closed convex rational $k$-dimensional polytopes $P_i\in\R^n$.
Its \textit{stable image}  under a rational vector map $p:\R^n\to\R^m$ is the union of those of $p(P_i)$ that are $k$-dimensional.
The \textit{stable intersection} of sets $P$ and $Q$ in $\R^n$ is the set of all points $x\in P\cap Q$, such that
$$\forall\varepsilon\, \exists\delta\, :\, v\in\R^n, |v|<\delta\, \Rightarrow\, \dist(x,(P+v)\cap Q)<\varepsilon.$$ This operation is denoted by $\wedge$,
is commutative, but not associative (different brackets in $\{x=0\}\wedge\{y=0\}\wedge\{y\leqslant|x|\}\subset\R^2$ lead to different answers), and
its result may be of unexpected dimension (like $\{z=|x|+|y|\}\wedge\{z=0\}=\{0\}\subset\R^3$). To avoid these issues, we should restrict our consideration to \textit{tropical complexes}. \par
A point $x$ of a $k$-dimensional polyhedral complex $P$ is said to be \textit{smooth}, if its transposed copy $P-x$ coincides with a vector subspace in a neighborhood of $0\in\R^n$. This vector
subspace is called the \textit{tangent space at}  $x$ and is denoted by $T_xP$.
\begin{defin} \label{deftropcompl} A closed polyhedral complex $P\subset\R^n$ is said to be \textit{tropical}, if it admits a positive 
locally constant non-zero function $w:\{$smooth points of $P\}\to\R$,
such that, for every rational subspace $L\subset\R^n$ of complementary dimension, the \textit{tropical intersection number}
of its transposed copy $L-x$ and $P$
$$\sum_{p\in P\cap(L-x)} w(p)\left|\frac{\Z^n}{(\Z^n\cap L) + (\Z^n\cap T_{p}P)}\right|$$ does not depend on $x$ (this sum makes sense for almost all $x\in\R^n$).
\end{defin}
Positivity of $w$ is not always included in the difinition, but will be crucial in what follows. We refer to the canonical papers \cite{fs} and \cite{st} or \cite{kaz} for the lemma below and for all subsequent well known facts about tropical complexes.
\begin{lemma} \label{lstab0} 1) Stable intersections, stable images and connected components of tropical complexes are tropical complexes. \newline 2) We have $\codim P\wedge Q=\codim P+\codim Q$ for tropical $P$ and $Q$, unless $P\wedge Q=\varnothing$.
\newline 3) Stable intersection is associative on tropical complexes.
\newline 4) The dual complex of a polytope is a tropical complex.
\end{lemma}
For a linear function $\gamma:\R^n\to\R$, we denote its support face of a polytope $\Delta$ by $\Delta^\gamma$ (it is the maximal face of $\Delta$ in the boundary of the half-space $\Delta+\{\gamma<0\}$), and say that $\gamma$ is an \textit{external normal covector} to the face $\Delta^\gamma$. 
Recall that the \textit{dual complex} of a polytope is the set of all external normal covectors to all of its positive-dimensional faces; in particular, the dual complex of a point is empty. Also recall that the \textit{mixed volume} is the unique real-valued multilinear (with respect to \textit{Minkowski summation } $A+B=\{a+b\, |\, a\in A ,\, b\in B\}$) function of $n$ convex polytopes in $\R^n$, whose value at $n$ copies of $A$ equals $n!\Vol(A)$. For any polytopes $A_1,\ldots,A_m$, we set the formal monomial $A_1\cdot\ldots\cdot A_m$ to be equal to 0 unless the affine span $L$ of $\sum_i A_i$ is $m$-dimensional; otherwise we set this monomial to be equal to the mixed volume of $A_1,\ldots,A_m$, induced by the integer volume form on $L$.
\par Let $\Delta_1,\ldots,\Delta_k$ be arbitrary polytopes in $\R^n$.
\begin{lemma}[\textbf{Equivalence lemma}] \label{lequiv1} The following conditions are equivalent:
\newline 1) For every $I\subset\{1,\ldots,k\}$, we have $\dim\sum_{i\in I}\Delta_i\geqslant|I|$.
\newline 2) There exist positive integers $a_1+\ldots+a_k=\dim\sum_i\Delta_i$, such that $\Delta_1^{a_1}\cdot\ldots\cdot\Delta_k^{a_k}\ne 0$.
\newline 3) The stable intersection of the dual complexes of $\Delta_1,\ldots,\Delta_k$ is non-empty.
\end{lemma}
\textsc{Proof.} For $n=k$, this equivalence is well known. In this case, the equivalence $1\Leftrightarrow 2$ was stated in \cite{kh77} (see e.g. \cite{dcg} for the proof), and the equivalence $2\Leftrightarrow 3$ follows from the fact that the mixed volume of the polytopes $\Delta_1,\ldots,\Delta_n$ equals the tropical intersection number of their dual complexes (see e.g. \cite{fs} or \cite{kaz}). In the general case, we may assume without loss of generality that $n=\dim\sum_i\Delta_i$, and have:\newline
$1\Rightarrow 2)\quad$ If the collection $\Delta_1,\ldots,\Delta_k$ satisfies the condition $(1)$, then the collection $$\Delta_1,\ldots,\Delta_k,\underbrace{\sum_i\Delta_i,\ldots,\sum_i\Delta_i}_{n-k}$$ also satisfiess the condition $(1)$. Since $1\Leftrightarrow 2$ for the latter collection, we have $\Delta_1\cdot\ldots\cdot\Delta_k\cdot(\Delta_1+\ldots+\Delta_k)^{n-k}>0$, or, removing the brackets, $\sum_{a_1+\ldots+a_k=n}C^n_{a_1,\ldots,a_k}\Delta_1^{a_1}\cdot\ldots\cdot\Delta_k^{a_k}>0$. Thus one of the terms in the latter sum is also strictly positive.\newline
$2\Rightarrow 3)\quad$ For the collection $\underbrace{\Delta_1,\ldots,\Delta_1}_{a_1},\ldots,\underbrace{\Delta_k,\ldots,\Delta_k}_{a_k}$, we have $2\Leftrightarrow 3$. Since this collection satisfies $(2)$, then the stable intersection of the dual fans of these $n$ polytopes is non-empty, thus the stable intersection of the dual fans of the subcollection $\Delta_1,\ldots,\Delta_k$ is also non-empty.\newline
$3\Rightarrow 1)\quad$ If we have $\dim\sum_{i\in I}\Delta_i<|I|$ for some $I$, then the polytopes $\Delta'_i=\Delta_i-a_i,\, a_i\in\Delta_i,\, i\in I$, are contained in an $I$-dimensional vector subspace $L\subset\R^n$, and the stable intersection of the dual fans of $\Delta'_i,\, i\in I$, has negative dimension in $L^*$ by Lemma \ref{lstab0}, i.e. is empty. Since the dual fans of $\Delta_i$ are the preimages of the dual fans of $\Delta'_i$ under the projection $(\R^n)^*\to L^*$, the stable intersection of them is also empty. $\quad\Box$
\begin{defin} Polytopes $A_1,\ldots,A_k$ are said to be \textit{linearly independent}, if they satisfy any of the three conditions of Lemma \ref{lequiv1}.
\end{defin}
We also need the following relative version of Equivalence lemma \ref{lequiv1}.
\begin{lemma} \label{lequiv_rel} Let $v:\Z^n\to\Z$ be a linear function such that the faces $\Delta_2^v,\ldots,\Delta_k^v$ are linearly independent, $\Delta_1^v$ is a point, and a (closed) polytope $\Delta$ is contained in $\Delta_1\setminus\Delta_1^v$. Then we have $$\sum_{a_1+\ldots+a_k=\dim\sum_i\Delta_i}\Delta_1^{a_1}\cdot\Delta_2^{a_2}\cdot\ldots\cdot\Delta_k^{a_k}>\sum_{a_1+\ldots+a_k=\dim\sum_i\Delta_i}\Delta^{a_1}\cdot\Delta_2^{a_2}\cdot\ldots\cdot\Delta_k^{a_k},$$ where the sum is taken over all decompositions of $\dim\sum_i\Delta_i$ into positive integers $a_i$.
\end{lemma}
\textsc{Proof.} \textbf{1.} With no loss in generality, assume that $\sum_i\Delta_i$ is not contained in an affine hyperplane, $\Delta_1^v=\{0\}$, and $\Delta$ contains all other vertices of $\Delta_1$. By a height-base formula for the mixed volume (see e.g. \cite{bernst}), we can rewrite $$\sum_{a_1+\ldots+a_k=\dim\sum_i\Delta_i}\Delta_1^{a_1}\cdot\Delta_2^{a_2}\cdot\ldots\cdot\Delta_k^{a_k}-\Delta^{a_1}\cdot\Delta_2^{a_2}\cdot\ldots\cdot\Delta_k^{a_k}=$$ $$=-\sum_\gamma \max_{x\in\Delta}\gamma(x)\cdot\sum_{a_1+\ldots+a_k=n}(\Delta^\gamma)^{a_1-1}(\Delta_2^\gamma)^{a_2}\ldots(\Delta_k^\gamma)^{a_k}, \eqno{(*)}$$ where $\gamma:\Z^n\to\Z$ runs over all primitive linear functions, whose restrictions to $\Delta$ are strictly negative. Choose primitive $\delta:\Z^n\to\Z$, whose restriction to $\Delta$ is strictly negative, so that $\Delta_2^\delta,\ldots,\Delta_k^\delta$ are linearly independent (there exists at least one such $\delta=v$), and $\dim(\Delta_2^\delta+\ldots+\Delta_k^\delta)$ is maximal possible. With no loss in generality, we also assume that the affine span of every $\Delta_i^\delta$ contains $0$.

\textbf{2.} If $\dim(\Delta_2^\delta+\ldots+\Delta_k^\delta)=n-1$, then, by Lemma \ref{lequiv1}, there exist positive $a_2,\ldots,a_k$, such that $a_2+\ldots+a_k=n-1$ and $(\Delta_2^\delta)^{a_2}\ldots(\Delta_k^\delta)^{a_k}>0$. Thus, in this case, the collection $\gamma=\delta,a_1=1,a_2,\ldots,a_k$ corresponds to the desired strictly positive term in $(*)$, and
it remains to consider the case $\dim(\Delta_2^\delta+\ldots+\Delta_k^\delta)<n-1$. In particular, 
let $N$ be the quotient space of $\R^n$ by the vector subspace, parallel to the affine span of $\Delta_2^\delta+\ldots+\Delta_k^\delta$,  then we can assume $\dim N>1$.

\textbf{3.}  Denote the projection $\R^n\to N$ by $q$ and the convex cone, generated by $q(\Delta_1)$, by $C$, then $q(\Delta_i^\delta)\subset C$ for every $i$. Otherwise there exists a linear function $\lambda$ such that $C^\lambda=\{0\}$, and the face $q(\Delta_2)^\lambda+\ldots+q(\Delta_k)^\lambda$ of positive dimension is adjacent to the vertex $q(\Delta_2^\delta+\ldots+\Delta_k^\delta)$, thus the linear function $\delta'=q^*(\lambda)$, as well as $\delta$, corresponds to linearly independent faces $\Delta_2^{\delta'},\ldots,\Delta_k^{\delta'}$ and is negative on $\Delta$. This would contradict the choice of $\delta$ on Step \textbf{1}, because $\dim(\Delta_2^{\delta'}+\ldots+\Delta_k^{\delta'})>\dim(\Delta_2^\delta+\ldots+\Delta_k^\delta)$.

\textbf{4.} The cone $C$ is not contained in a hyperplane (otherwise every $q(\Delta_i)$ is contained in the same hyperplane $H$ by Step \textbf{3}, and $\Delta_1+\ldots+\Delta_k$ is contained in the hyperplane $q^{-1}(H)$, which contradicts the assumption). Thus, there exists a linear function $\lambda:N\to\Z$ such that $\dim q(\Delta)^\lambda=\dim N-1$, which is positive by Step \textbf{2}, and $C^\lambda=\{0\}$. Taking $\gamma=q^*\lambda$, we conclude that $\Delta_2^\gamma,\ldots,\Delta_k^\gamma$ are linearly independent, $\dim(\Delta^\gamma+\Delta_2^\gamma+\ldots+\Delta_k^\gamma)=n-1$, and $\gamma$ is negative on $\Delta$. By Lemma \ref{lequiv1}, there exist positive $a_1,\ldots,a_k$, such that $a_1+\ldots+a_k=n$ and
$(\Delta^\gamma)^{a_1-1}(\Delta_2^\gamma)^{a_2}\ldots(\Delta_k^\gamma)^{a_k}>0$. Thus, the collection $\gamma,a_1,\ldots,a_k$ corresponds to the desired strictly positive term in $(*)$.
$\quad\Box$\par

The following fact can be seen as a tropical version of purity results for bifurcation sets over $\C$, and will eventually allow to prove these results themselves.
\begin{utver} \label{dtrop} Let $P$ be a polyhedron in $\R^n$, open in its affine span, and let $T$ be a tropical complex in $\R^n$ such that $\dim P+\dim T=n+k$. Then \newline 1) The stable intersection $S=T\wedge P$ is $k$-dimensional or empty. \newline
2) The intersection of the closure of $S$ with the relative boundary $\partial P$ is $(k-1)$-dimensional or empty.
\newline 3) It is empty if and only if every connected component of $S$ is contained in an affine subspace that is contained in $P$.
\end{utver}
Both statements remain valid, if we define $S$ as the conventional intersection $T\cap P$, and claim the dimension of the intersections in $(1)$ and $(2)$ to be greater or equal than what we have in the stable case. 
The proof of this refinement follows the same lines as the proof of Proposition \ref{dtrop}. 

Before proving Proposition \ref{dtrop}, we explain in what sense it is the tropical version of purity results for bifurcation sets over $\C$. Let $T$ be a $p$-dimensional tropical complex in $\R^q\times\R^p$. A point $x\in T$ is said to be {\it regular} for the projection $T\to\R^p$, if a generic fiber of this projection has at most one point in a small neighborhood of $x$. The {\it tropical Jelonek set} of the projection $T\to\R^p$ is the set of images of all points $x\in T$ that are not regular for this projection. A $p$-dimensional tropical complex $T\subset\R^{p+q}$ is said to be {\it regular}, if every its point $x$ admints a projection $\R^{p+q}\to\R^p$, for which the point $x$ is regular.
\begin{conj}[Tropical Jelonek theorem] \label{tropjel} If a $p$-dimensional tropical complex is regular,
then the tropical Jelonek set of every its projection to $\R^p$ is a polyhedral (not necessarily tropical) complex of pure codimension 1.
\end{conj}
Let $P$ be a convex polyhedron in $\R^n$, represent it as $\{x\, |\, h(x)=c\}$, where $c\in\R$, and $h:\R^n\to\R$ is a continuous piecewise linear function, whose restriction to every ray from the origin is linear. Let $T\in\R^1\times\R^n$ be the corner locus of the function $\max(|y|,h-c)$ on $\R^1\times\R^n$, where $y$ is the standard coordinate on $\R^1$. Then the tropical Jelonek theorem for the projection $T\to\R^n$ is exactly Proposition \ref{dtrop}.

\begin{lemma}\label{l27} 1) If the image of a tropical complex $T$ under a linear map $\R^n\to\R^k$ is contained in a polytope $\Delta$ that does not contain a line, then this image consists of finitely many points. In particular, a positive-dimensional tropical complex is unbounded.
\newline 2) If, for every smooth point $x$ of a tropical complex $T\subset\R^n$, the tangent space $T_x T$ contains a vector subspace $L$, then $T$ is the preimage of a tropical complex in $\R^n/L$.
\end{lemma} 
\textsc{Proof of Part 1.} Assume that  the image is of positive dimension, then one can find a polytope $T_i\subset T$, such that its image in $\R^k$ is at least 1-dimensional, then there exists an affine subspace $L\subset\R^n$ of complementary dimension, intersecting $T_i$ transversally, such that the intersection of its image in $\R^k$ with $\Delta$ is bounded. Thus the tropical intersection number of $L$ with $T$ (see Definition \ref{deftropcompl}) is positive, while, for suitable $x\in\R^n$, the image of $x+L$ in $\R^k$ does not intersect $\Delta$, and the tropical intersection number of $x+L$ with $T$ is 0.


\textsc{Proof of Part 2.} Assume that $T$ is not the preimage of a subset in $\R^n/L$: then there exists a point $x\in T$, such that a small neighborhood of $x$ in $L+x$ is not contained in $T$. With no loss in generality, we can assume that $x=0$, and consider the union of polyhedral cones $T'$ that coincides with $T$ in a small neighborhood of the origin. Since the image of $T'$ in $\R^n/L$ is of dimension $\dim T-\dim L$, there exists a vector subspace $L'$ of complementary dimension, whose image in $\R^n/L$ intersects the image of $T'$ at one point $0$. Since $0\in T'\cap L'$, the tropical intersection number of $L'$ and $T'$ is positive,
however, since $L$ is not contained in $T'$, we can choose $x'\in L\setminus T'$ and notice that $T'\cap(L'+x')=\varnothing$, thus the the intersection number of $T'$ and $L'+x'$ is 0, thus $T'$ is not tropical, thus $T$ is not tropical. 

In the same way one can check that a function $w:\{$smooth points of $T\}\to\R$, defining a tropical structure on $T$ (see Definition \ref{deftropcompl}), is constant along every affine space, parallel to $L$, and thus induces a function on the smooth part of the image of $T$ in $\R^n/L$, which defines a tropical structure on the image. $\quad\Box$

\textsc{Proof of Proposition \ref{dtrop}.} Let $L$ be the affine span of $P$, then $T\wedge P = (T\wedge L)\wedge P = (T\wedge L)\cap P$. Thus, Part 1 follows by Lemma \ref{lstab0}, and, replacing $T$ with $T\wedge L$ and $\R^n$ with $L$ ,
we may assume that $P$ is an open $n$-dimensional polytope in $\R^n$. 
We also assume by induction that the desired statement is already proved for polytopes of dimension smaller than $n$. 

\textsc{Proof of Part 2.} Assume the contrary to the statement of Part 2: there exists an open polytope $P$, a tropical complex $T$ in $\R^n$ and a point $x\in \overline{T\cap P}\cap\partial P$, such that the set $\overline{T\cap P}\cap\partial P$ is of dimension less than $\dim T-1$ in a small neighborhood of $x$ (in particular, $\dim T>1$). With no loss in generality, we can also assume that $x=0$, consider
the polyhedral cone $P'$ that coincides with $P$ in a small neighborhood of the origin, and the tropical complex $T'$ that consists of finitely many polyhedral cones $T_i$ and coincides with $T$ in a small neighborhood of the origin. Then, by our assumptions, the set $\overline{T'\cap P'}\cap\partial P'$ is non-empty of dimension smaller than $\dim T-1$. 

We now prove that $P'$ contains a line. If it does not, then there exists an affine hyperplane $L''$, intersecting it at a bounded polytope $P''$ of dimension $n-1$. Denote the (stable) intersection $T'\wedge L''=T'\cap L''$ by $T''$,
then $\overline{T''\cap P''}\cap\partial P''=L''\cap \overline{T'\cap P'}\cap\partial P'$ is of dimension smaller than $k-2$. On the other hand, $\overline{T''\cap P''}\cap\partial P''$ is not empty, otherwise the connected component of $T'\wedge L''$ in $P''$ is a non-empty bounded positive-dimensional tropical complex, which is impossible by Lemma \ref{l27}.1. Thus, the smaller-dimensional polytope $P''$ and the tropical complex $T''$  contradict the inductive assumption. Thus, the maximal vector subspace $L'$, contained in $T'$, is positive-dimensional. 

We now prove that $L'$ is parallel to every $T_i$. If it does not, then there exists a line $l\subset L'$ that is not parallel to some $T_i$, then the stable image $T''$ of $T'$ under the projection along $l$ is non-empty. Denote the projection of $P'$ along $l$ by $P''$, then the set $\overline{T''\cap P''}\cap\partial P''$ contains $0$, is contained in the projection of the set $\overline{T'\cap P'}\cap\partial P'$,  and its dimension is also smaller that $k-1$, thus the smaller-dimensional polytope $P''$ and the tropical complex $T''$ contradict the inductive assumption.

We now prove that $L'$ is not parallel to some $T_i$. Otherwise, by Lemma \ref{l27}.2, $T'=T''\times L'$ and $P'=P''\times L'$ for a polytope $P''$ and a tropical complex $T''$ in the quotient space $\R^n/L'$. Since $\overline{T'\cap P'}\cap\partial P'=(\overline{T''\cap P''}\cap\partial P'')\times L''$,
the dimension of $\overline{T''\cap P''}\cap\partial P''$ is less than $\dim T''-1$, thus the smaller-dimensional polytope $P''$ and the tropical complex $T''$ contradict the inductive assumption.

Assuming the contrary to the statement of Part 2, we have concluded that $L'$ is parallel to every $T_i$, but not parallel to one of them, which is a contradiction.

\textsc{Proof of Part 3.} Represent the polytope $P$ as the preimage of a polytope $\Delta\subset\R^k$ under a linear map $\R^n\to\R^k$, so that $\Delta$ does not contain a line, and apply Lemma \ref{l27}.1.$\quad\Box$

\subsection{Essential tuples and mixed fiber polytopes} \label{Smixfib}

Many polytope-related notions of convex geometry admit a mixed version (e.g. mixed volume, mixed fiber polytope), related to a tuple of polytopes. We aim at inventing a mixed version of the notion of face, and need the following vector notation for tuples of polytopes.

Let $A$ be a $K$-tuple of polytopes in $\R^n$ (i.e. a map from a finite
set $K\in\Z$ to the set of polytopes in $\R^n$). For every
$J\subset K$, denote the restriction of $A$ to $J$ by $AJ$ (it is
called a \textit{subtuple} of $A$), the cardinality of $K$ by $|A|$,
the Minkowski sum of $A(k),\, k\in K$, by $\sum A$ (which is
$\{0\}$ for $K=\varnothing$), the difference $\dim(\sum A) - |A|$
by $\dim A$, the minimum of $\dim AJ$
over all $J\subset K$ by $\min\dim A\leqslant 0$, the mixed volume or mixed fiber
polytope of the polytopes $A(k),\, k\in K$, by $\MV(A)$ or $\MP(A)$ respectively,
the tuple of the images
of $A(k),\, k\in K\setminus J$, under the projection along the affine span
of $\sum AJ$ by $A/J$. For a tuple of numbers $a:K\to\Z$, let $A^a$
be the tuple that includes $a(k)$ copies of the polytope $A(k)$
for every $k\in K$ (more precisely, choose a map $\kappa:K'\to K$,
such that the preimage of every $k\in K$ consists of $a(k)$ elements,
and set $A^a$ to be the $K'$-tuple $A\circ\kappa$).
\begin{lemma} \label{mvprod} For a tuple $A$ of linearly independent polytopes, such that $\dim A=\dim AI=0$, we have $\MV(A)=\MV(AI)\MV(A/I)$.
\end{lemma} See e.g. \cite{dcg} for the proof.
\begin{defin}[\cite{sturmf}]\label{defess} The tuple $A$ is said to be \textit{essential,} if $\dim AJ=\min\dim A$ only for $J=K$. \end{defin}
\begin{lemma}[\cite{sturmf}]\label{ldefess} 1) There exists a unique minimal $J\subset K$ with $\dim AJ=\min\dim A$; \newline 2) There exists a unique maximal $J\subset K$ with $AJ$ essential. \newline 3) These two subsets coincide.\end{lemma}
\textsc{Proof.} Note that $$\dim AI+\dim AJ\geqslant \dim A(I\cup J)+\dim A(I\cap J)\eqno{(*)}$$ (see \cite{sturmf}). If $\dim AI=\dim AJ=\min\dim A$ then $\dim A(I\cap J)=\min\dim A$ as well by $(*)$, then the intersection of all $I$ with $\dim AI=\min\dim A$ has the same dimension, and $(1)$ follows. We denote this intersection by $I_0$ and note that $AI_0$ is essential. 

If 
$AJ$ is essential, and $J\not\subset I$, then $\dim A(I\cap J)>\dim AJ$, and, summing up this equality and $(*)$, we get $\dim AI>\dim A(I\cup J)\geqslant\min\dim A$, which implies that $I\ne I_0$. Thus, $I_0$ contains any other $J$ with $AJ$ essential, which proves $(2)$ and $(3).\quad\Box$

The first application of this lemma is a version of Equivalence lemma \ref{lequiv1} for mixed fiber polytopes, whose definition we recall now. 
Let $P_L$ be the semigroup of polytopes in a vector space $L$ (with respect to Minkowski summation).
For a polytope $H\subset\R^n\oplus\R^m$, the \textit{fiber polytope} $\int H$ is the set of points of the form
$\int_{\R^n}\gamma(x)dx\subset\R^m$, where $\gamma:\R^n\to\R^m$ runs over all continuous sections of the projection $H\to\R^n$ (\cite{bs}).
\begin{defin}[\cite{mcm}, \cite{ekh}]\label{defmp} The (unique) symmetric multilinear mapping $$\MP:\underbrace{P_{\R^n\oplus\R^m}\times\ldots\times P_{\R^n\oplus\R^m}}_{n+1}\to P_{\R^m},$$ such that $\MP(H,\ldots,H)=(n+1)!\int H$ for every $H$, is called the \textit{mixed fiber polytope}.
\end{defin}
The multiplier $(n+1)!$ assures that the mixed fiber polytope of integer polytopes is integer.
\begin{lemma}[\cite{ekh}]\label{lekh} The polytope $X$ satisfies the equality
$$\MV(X,B_1,\ldots,B_{m-1})=\MV(H_0,\ldots,H_n,B_1,\ldots,B_{m-1})$$
for all tuples of polytopes $B_1,\ldots,B_{m-1}\subset\R^m$, iff $X=\MP(H_0,\ldots,H_n)+a$ for some $a\in\R^m$.
\end{lemma}
In the left hand side of this equation $\MV$ denotes the mixed volume of $m$ polytopes in $\R^m$,
and in the right hand side -- the same for $m+n$ polytopes in $\R^n\oplus\R^m\supset\{0\}\times\R^m\simeq\R^m$.

Let $H$ be a $\{0,\ldots,k\}$-tuple of polytopes in $\R^n\oplus\R^m$, and $A$ be its projection to $\R^n$ (i.e. $A(i)\subset\R^n$ is the projection of $H(i)$ for every $i$).
Assume that $A$ is relevant (Definition \ref{defrel}). 
\begin{theor} \label{trivmp} 
1) Assume that $A$ is linearly dependent, and $AI$ is its maximal essential subtuple.
The sum $$\sum_{a_0+\ldots+a_k=n+1}\MP(H^a)\eqno{(*)}$$ over
all decompositions of $n+1$ into a sum of positive integers is a point, iff the projection of $\sum HI$ to $\R^n$ is injective.
\newline 2) Assume that $A$ is linearly independent.
The sum $(*)$ is a point, iff the projection of $\sum H$ to $\R^n$ is injective.
\end{theor}
\textsc{Proof of the implication $\Leftarrow$.} In both cases, under given assumptions, the mixed volume of the polytopes
$$\underbrace{H(0),\ldots,H(0)}_{a_0},\ldots,\underbrace{H(k),\ldots,H(k)}_{a_k},B_1,\ldots,B_{m-1}$$
is zero by Lemma \ref{lequiv1} for every decomposition of $n+1$ into the sum of positive integers $a_0+\ldots+a_k$, and for all polytopes $B_0,\ldots,B_{m-1}$ in $\R^m$.
By Lemma \ref{lekh}, this implies that the mixed volume of $B_1,\ldots,B_{m-1}$ and the mixed fiber polytope $\MP(H^a)$
is zero, thus this polytope is a point.

\textsc{Proof of $\Rightarrow$.} 
Consider a generic set of segments $B_1,B_2,\ldots, B_{m-1}\subset \R^m$. Let $L$ be the affine span of $\sum B_i$.

Consider the $\{0,1,\ldots, k+m-1\}$-tuple $R$ such that $R(i)=H(i), i\leq k$ and $R(i)=B_{i-k}, i>k$. Consider its subtuple $B:=R\{k+1,\ldots, k+m-1\}$.

If the projection of $\sum H$ to $\R^n$ is injective, there is nothing to prove. Otherwise, chose a one-dimensional subspace $K\subset \R^m$ that is parallel to the affine span of a generic fiber. Since $L$ is generic, we have $\R^m = L\oplus K$, and the space $(\R^m \oplus \R^n)/L$ is naturally isomorphic to $K\oplus \R^n$. Thus we have $\dim \sum R/B = n+1$, and therefore $\dim \sum R = \dim \sum R/B + \dim \sum B = n+m$. We check the conditions of Lemma \ref{lequiv1} as applied to $R$. Since the sum (*) is a point, we have $\sum_{a_0+a_1+\ldots+a_k=n+1} H_0^{a_0}\ldots H_n^{a_n}B_1\ldots B_{m-1} = 0.$ As the polyhedra $B_i$ are one-dimensional, it follows that
$$
\sum_{a_0+a_1+\ldots+a_{k+m-1}=n+m} \MV(R^a) = 0.
$$
According to Lemma \ref{lequiv1} there is a subtuple $RJ$ such that $\dim RJ <0$. 

We can assume that $B$ is a subtuple of $RJ$, otherwise we replace $J$ with $J\cup \{k+1,\ldots,k+m-1\}$, and $\dim RJ <0$ remains valid. Let $J'=J\setminus \{k+1,\ldots,k+m-1\}$. Let $p(HJ')$ denotes the projection of $HJ'$ to $K\oplus \R^n$. We have $p(HJ') = RJ/B$ and $\dim p(HJ') = \dim RJ/B = \dim RJ <0$. 

Assume that the projection of $\sum HJ'$ to $\R^n$ is not injective. Since $L$ is generic, a fiber of the projection is transversal to $L\oplus \R^n$. Therefore we have $\dim p(HJ') = \dim AJ' + 1$ and thus $\dim AJ' < -1$, which contradicts the condition that $A$ is relevant. Thus, the projection of $\sum HJ'$ to $\R^n$ is injective and $\dim AJ' = \dim p(HJ') =-1$. Thus, $J'$ contains $I$ such that $AI$ is the maximal essential subtuple of $A$, and the projection of $\sum HI$ to $\R^n$ is also injective. $\quad\Box$

\subsection{Essential faces and their adjacency} \label{Sessfaces}
\begin{defin} A tuple of polytopes $B$ is called a \textit{face} of the tuple $A$, if $B(i)$ is a face of $A(i)$ for every $i\in K$, and $\sum B$ is a face of $\sum A$. It is written as $B\prec A$ and said to be \textit{trivial}, if $B=A$.
\end{defin}
\begin{defin} \label{dessf} A tuple $E$ is called an
\textit{essential facing} of the tuple $A$, if it is the maximal
essential subtuple of a face $B\prec A$. It is said to be \textit{trivial},
if it is a subtuple of $A$. An essential facing $E$ is said to be \textit{adjacent} to an essential facing $E'$, if they are maximal essential subtuples of faces $B\prec B'\prec A$ respectively.
\end{defin}
\begin{theor}\label{thess_main} Every essential facing $E$ is adjacent to another essential facing $E'$ such that $\dim E'=\dim E+1$.
\end{theor}
Thus, essential facings of a tuple of polytopes behave similarly to faces of one polytope in terms of adjacency and dimension. On the other hand, adjacency of essential facings is not transitive (see the following example), which is the primary source of complications in this section.
Before proving the statement, we give counterexamples to a number of tempting, but misleading ways to simplify it.
Denote the standard simplex and the standard cube in $\R^n$ by $S_n$ and $C_n$, and denote the convex hull of the points $(0,0,0),(0,1,2),(0,2,1)$ by $D$.
\begin{exa} In the statement of Theorem, the essential facing $E$ is not necessarily a face of a subtuple of $E'$: consider $A(1)=C_2, A(2)=A(3)=A(1)^{(0,-1)}$ and $E=A^{(0,-1)}$, then $E'=A\{2,3\}$. The tuple $A$ is shown in black on the left of the picture below, with $E$ in blue and $E'$ in red.

In the statement of Theorem, a subtuple of the essential facing $E$ is not necessary a face of $E'$: consider $A(1)=C_2, A(2)=A(3)=S_2$ and $E=A\{2,3\}^{(0,1)}$, then $E'=A$. The tuple $A$ is shown in black on the right of the picture below, with $E$ in blue.

\noindent\includegraphics[width=\textwidth]{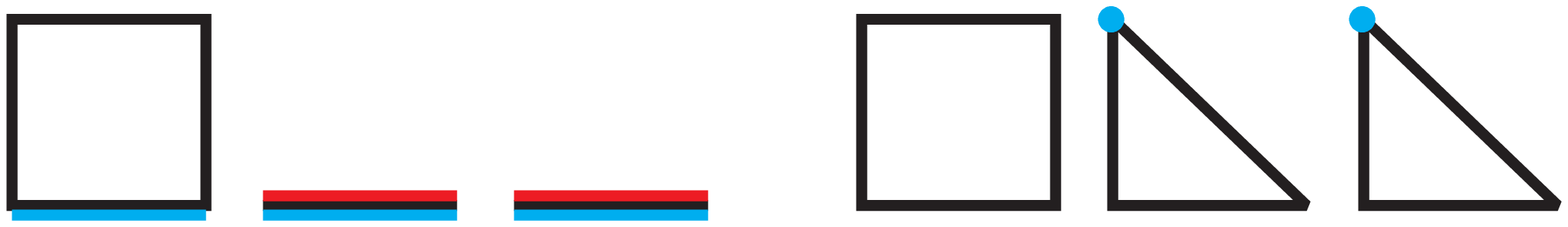}

If $V, E$ and $E'$ are essential facings of $A$ such that $V$ is adjacent to $E$ and $E$ is adjacent to $E'$, then $V$ may be not adjacent to $E'$: consider $A(1)=C_3, A(2)=A(3)=A(4)=D$, $V=A^{(-1,-1,-1)}$, $E=A\{2,3,4\}^{(0,-1,-1)}$ and $E'=A\{2,3,4\}^{(0,-2,1)}$. The tuple $A$ is shown in black on the picture below,
with $V$ shown by blue dots, $E$ shown by red rounds, and $E'$ shown by green segments. Some possible ways to represent $E$ or $E'$ as a maximal essential subtuple of a face $B$ or $B'\prec A$ are shown in red and green dotted lines respectively.
\begin{center}
\noindent\includegraphics[width=0.8\textwidth]{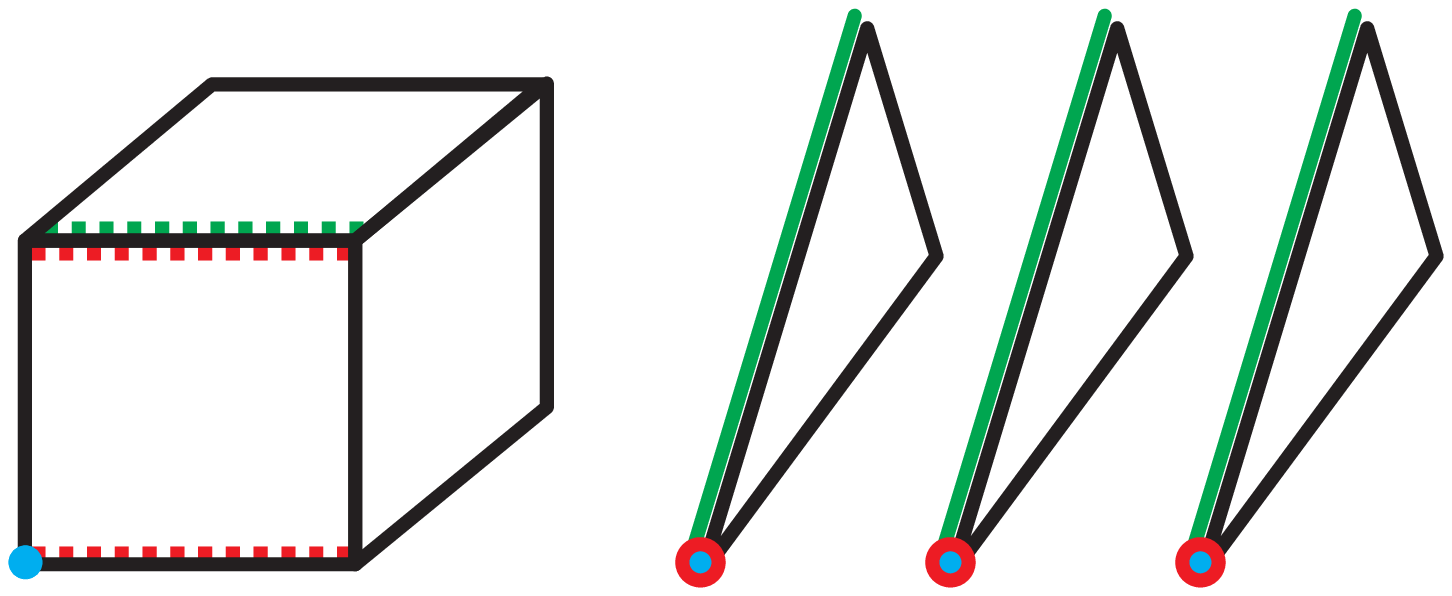}\end{center}
\end{exa}
When working out these examples, it is useful to reformulate the definition of essential facings and their adjacency as follows (although we do not use this statement in the sequel and omit details of its proof).
\begin{utver} \label{altdefessface} 1) Every essential facing $E$ of a tuple $A$ is the maximal essential subtuple of a face $B\prec A$ such that $\dim E=\dim B$. \newline 2) Every pair of adjacent facings $E$ and $E'$ are maximal essential subtuples of faces $B\prec B'\prec A$ such that $\dim E=\dim B$ and $\dim E'=\dim B'$.
\end{utver}
Part 1 is Lemma \ref{thess00} below, and Part 2 can be proved by induction on $\dim E'-\dim E$, using Theorem \ref{thess1} below. 

However, given $B$ that satisfies the statement of Part 1, it maybe impossible to represent it as a face of $B'$, satisfying the statement of Part 2, even if $\dim E'=\dim E+1$. 
For instance, in the notation of the last example above, the tuple $B=A^{(0,-1,-1)}$ proves the statement of Proposition \ref{altdefessface}(1) for $E$ as introduced in the example, but does not admit $B'$ such that $B$ and $B'$ would prove Proposition \ref{altdefessface}(2) for $E$ and $E'$ as introduced in the example (it is proved by $B=A^{(0,-10,1)}$ and $B'=A^{(0,-2,1)}$ instead).

We prove the following stronger version of Theorem \ref{thess_main}.
\begin{theor} \label{thess1}
If $E$ is a nontrivial essential facing of $A$,
then it is also a nontrivial essential facing of a face $B\prec A$,
such that $\dim B=\min\dim B=1+\dim E$.
\end{theor}

The rest of the section is devoted to the proof of this theorem. We need the following relation between essential tuples and stable intersections:
\begin{utver} \label{stabess} For a decomposition $K=I\sqcup H$, let $L$ be the orthogonal complement to the sum $\sum AI$,
and let $C$ be the stable intersection of the dual complexes of the polytopes $AH$.
Then \newline 1) If $AI$ is the maximal essential subtuple of $A$,
then the stable intersection of $L$ and $C$ is non-empty. \newline 2) If the stable intersection of $L$ and $C$ is non-empty,
then the maximal essential subtuple of $A$ is a subtuple of $AI$. \end{utver}
Note that we could not reverse any of these two statements.
\par
\textsc{Proof of Part 1.} Choose arbitrary polytopes $B_1,\ldots,B_{\codim L}$, whose orthogonal complement equals $L$, then the tuple $B_1,\ldots,B_{\codim L},AH$ satisfies Condition (1) of Lemma \ref{lequiv1}. Then, by  this lemma, the stable intersection of the dual complexes of $B_1,\ldots,B_{\codim L},AH$ is non-empty, while the stable intersection of the dual complexes of $B_1,\ldots,B_{\codim L}$ equals $L$. \par

\textsc{Proof of Part 2.} Assume, on the contrary, that there exists an essential subtuple $AJ,\, J\not\subset I$.
Denote the orthogonal complement of $A(I\cap J)$ by $M$ and choose arbitrary polytopes $B_1,\ldots,B_{\codim M}$, whose orthogonal complement equals $M$. Since $AJ$ is essential, we have $$\dim\bigl(\Sigma A(H\cap J)+\Sigma_i B_i\bigr)=\dim AJ<|I|-|I\cap J|+\dim\Sigma A(I\cap J)=|H\cap J|+\dim\Sigma_i B_i,$$ i.e. the tuple $B_1,\ldots,B_{\codim M},A(H\cap J)$ does not satisfy Condition (1) of Lemma \ref{lequiv1}, thus the stable intersection of its dual complexes is empty. Since the stable intersection of the dual complexes of $B_1,\ldots,B_{\codim M}$ equals $M$, then it contains $L$, thus the stable intersection of $L$ and the dual complexes of $AH$ is empty as well. $\quad\Box$ \par

The following is well-known: \begin{lemma}\label{ltrunc0} The dual complex of $\Delta^\gamma$ is defined by the following conditions: it consists of cones with vertices at $\gamma$ and coincides with the dual complex of $\Delta$ in a neighborhood of $\gamma$. 
\end{lemma} For the tuple of polytopes $A$, denote its face $B$ with $B(i)=A(i)^\gamma,\, i\in K$, by $A^{\gamma}$. Note that every face can be represented in this form for a suitable $\gamma$. We now prove the following two weak versions of Theorem \ref{thess1}.
\begin{lemma} \label{thess00} If $E$ is the maximal essential subtuple of $A$, then it is also the maximal essential subtuple of a certain $B\prec A$ such that $\dim B=\dim E$.
\end{lemma}
\textsc{Proof.} Applying Proposition \ref{stabess}(1) with $AI=E$, we find that $L\wedge C$ is non-empty (in the notation of the proposition). Assuming that $\dim A>\dim E$ (otherwise we could set $B=A$), Lemma \ref{lstab0}(2) implies that $\dim L\wedge C>0$; in particular, $L\wedge C$ contains a non-zero covector $\gamma$. By Lemma \ref{ltrunc0}, the stable intersection of the dual complexes of $A(i)^\gamma,\, i\notin I,$ and $L$ is also non-empty (namely, it also contains $\gamma$). Thus, applying Proposition \ref{stabess}(2) to $A^\gamma$ and its subtuple $E$, the maximal essential collection of $A^\gamma$ is a subtuple of $E$. Since $E$ is essential itself, it is exactly the maximal essential collection of $A^\gamma$. We have proved the following: if $E$ is the maximal essential subtuple of $A$, such that $\dim A>\dim E$, then it is also  the maximal essential subtuple of a non-trivial face $B\prec A$. Iterating this, we get the desired statement.
$\quad\Box$ \par
\begin{lemma} \label{thess0}
If $E$ is a nontrivial essential facing of $A$,
then it is also a nontrivial essential facing of a nontrivial face $B\prec A$,
under an obviously necessary assumption $\dim A>1+\dim E$.
\end{lemma}
\textsc{Proof.} By Lemma \ref{thess00}, we can choose a face $D\prec A$, such that $E=DI$ is its maximal essential subtuple, and $\dim D=\dim E$. Applying Proposition \ref{stabess}(1) to $D$ and its maximal essential subtuple $E$, we get $C\wedge L\ne\varnothing$ (in the notation of the proposition).

We now prove the same statement for the stable intersection $C'$ of the dual complexes of $A(i),\, i\notin I,$ and the (relatively open) cone $L'\subset L$ of external normal covectors to the face $\sum E$ of the polytope $\sum AI$: the set $L'\wedge C'$ is non-empty as well. Assume the contrary: $L'\wedge C'=\varnothing$. Choose $\gamma\in L'$ so that $A^\gamma=D$,
then, by Lemma \ref{ltrunc0}, the complex $C$ is uniquely defined by the fact that it 
coincides with $C'$ in a small neighborhood of $\gamma$. Since $L$ 
also coincides with $L'$ in a small neighborhood of $\gamma$, we have $C\wedge L=L'\wedge C'=\varnothing$ in a small neighborhood of $\gamma$, and this implies $C\wedge L=\varnothing$.


By Lemma \ref{lstab0}(2) and the assumption $\dim A>1+\dim E$, we have $\dim L'\wedge C'>1$.
Applying Proposition \ref{dtrop}(2) and (3) to the tropical complex $C'$ and the polyhedron $L'$, we conclude that the closure of $L'\wedge C'$ intersects the (relative) boundary of $C'$ by a non-empty positive-dimensional set. In particular, this intersection contains a non-zero covector $\gamma$ that has a nearby non-zero covector $\delta\in L'\wedge C'$. We will prove that $A^{\gamma}$ is the desired face $B$ as follows. By Lemma \ref{ltrunc0}, the stable intersection of $L$ and the dual complexes of $A(i)^\delta,\, i\notin I,$ is non-empty (namely, it contains $\delta$, because $L'\wedge C'$ does). Thus, applying Proposition \ref{stabess}(2) to $A^\delta$ and its subtuple $E$, the maximal essential collection of $A^\delta$ is a subtuple of $E$. Since $E$ is essential itself, it is exactly the maximal essential collection of $A^\delta\prec A^\gamma$. $\quad\Box$\par
\begin{lemma}\label{l21} If $E$ is a trivial essential facing of $A$, then $\dim E=\min\dim A$.\newline If $E$ is a non-trivial essential facing of $A$, then $\dim E<\min\dim A$.
\end{lemma} \textsc{Proof.} The first statement follows by definition. To prove the second one, let $E$ be the maximal essential subtuple $BJ$ for a face $B\prec A$. If $J\subset I$, then $\dim E\leq\dim BI<\dim AI$, with the second inequality strict because $E$ is non-trivial. Otherwise, we have $\dim E<\dim BI\leq\dim AI$, with the first inequality strict by Lemma \ref{ldefess}(1). $\quad\Box$\par
\textsc{Proof of Theorem \ref{thess1}} proceeds by induction on $\dim A$. If $\dim A>1+\dim E$, then substitute $A$ with $B$ given by Lemma \ref{thess0}, and apply the inductive hypothesis. If $\dim A=1+\dim E$, then it is also equal to $\min\dim A$ by Lemma \ref{l21}, and we can set $B=A$. $\quad\Box$

\subsection{Lattice subsets and resultants} \label{latticesets}
We now formulate the results of the previous subsections in the form that we need later, with finite subsets of $\Z^n$ instead of polytopes. In particular, we apply all the polytope terminology to finite subsets of $\Z^n$ by considering convex hulls of these subsets.
A \textit{face of a finite subset} $H\subset\Z^n$ is the intersection of $H$ with a face of its convex hull. The dimension $\dim H$ of $H\subset\Z^n$ is the dimension of its convex hull.

As before, a $K$-\textit{tuple} of finite sets is a map $A$ from a finite set $K$ to the set of finite subsets of $\Z^n$, its $I$-\textit{subtuple} for $I\subset K$ is denoted by $AI$, the Minkowski sum of its elements $\{\sum_{k\in K}a_k\,|\,a_k\in A_k\}$ is denoted by $\sum A$ (which is $\{0\}$ for $K=\varnothing$), 
and the difference $\dim(\sum A)-|K|$ by $\dim A$. A $K$-tuple $\Gamma$ is a \textit{face} of a $K$-tuple $A$, if $\Gamma(k)$ is a face of $A(k)$ for every $k\in K$, and $\sum\Gamma$ is a face of $\sum A$.

\begin{defin} \label{defimport}
A face $\Gamma\prec AI$ is said to be \textit{important}, if there exists a face $\Gamma'\prec A$ such that $\Gamma=\Gamma'I$ and $\dim\Gamma\leqslant\dim\Gamma'J$ for every $J\supset I$. It is said to be \textit{essential}, if, in addition, $\dim\Gamma<\dim\Gamma'J$ for every $J\subsetneq I$ (i.e. if the convex hulls of $\Gamma$ form an essential facing of the convex hulls of $A$).
\end{defin}

If $\Gamma$ is a face of $AI$, then we denote the projection of $\Z^n$ along the affine span of $\sum\Gamma$ by $p$, the convex hulls of $pA(i),\, i\notin I$, by $\Psi_1,\ldots,\Psi_m$, the convex hulls of $p(\sum AI)$ and $p(\sum AI\setminus\sum\Gamma)$ by $\Psi_0$ and $\Psi$.
\begin{defin} \label{mixmilnum} The \textit{Milnor number} of the tuple $A$ at the face $\Gamma$ is $$c_A^\Gamma=\sum_{a_0,\ldots,a_m}\Psi_0^{a_0}\cdot\Psi_1^{a_1}\cdot\ldots\cdot\Psi_m^{a_m}-\Psi^{a_0}\cdot\Psi_1^{a_1}\cdot\ldots\cdot\Psi_m^{a_m},$$
where $(a_0,\ldots,a_m)$ runs over all collections of positive integers that sum up to $\dim\sum A-\dim\sum\Gamma$, and the mixed volume of polytopes $B_1,\ldots,B_k$ in a $k$-dimensional vector space is denoted by the monomial $B_1\cdot\ldots\cdot B_k$.\end{defin}
This number does not depend on all the tuple $A$, but only on its local structure near $\Gamma$. More precisely, let $\psi_0$ be the vertex $p(\sum\Gamma)$ of $\Psi_0$, and let $$C=\{ \lambda (\psi-\psi_0) \, |\, \lambda>0,\, \psi\in\Psi_0\}$$ be the cone that looks near the origin like $\Psi_0$ near its vertex $\psi_0$, and denote the pair of (unbounded) polyhedra $(\Psi_0+C,\Psi+C)$ by $\Phi_0$ and $(\Psi_i+C,\Psi_i+C)$ by $\Phi_i$ for $i>0$.
\begin{utver} The number $c_A^\Gamma$ equals $$\sum_{a_0,\ldots,a_m}\Phi_0^{a_0}\cdot\Phi_1^{a_1}\cdot\ldots\cdot\Phi_m^{a_m},$$
where $(a_0,\ldots,a_m)$ runs over all collections of positive integers that sum up to $\dim\sum A-\dim\sum\Gamma$, and the mixed volume of pairs of polytopes $P_1,\ldots,P_k$ in a $k$-dimensional vector space is denoted by the monomial $P_1\cdot\ldots\cdot P_k$ (see \cite{me06} or \cite{detandnewt} for its definition).
\end{utver}
See Assertion 1.5 in \cite{detandnewt} for the proof. 
We do not need this representation of Milnor numbers in what follows, we mention it only because it simplifies the computation of $c^\Gamma_A$ in many special cases (see \cite{detandnewt} for properties and examples of computation of the mixed volume of pairs of polyhedra). This number and its relation to discriminants is also discussed in detail in \cite{st11} (where it is called the {\it mixed multiplicity}), and an illustrative example for $n=2,\, k=1,\, \Gamma=($point, point$)$ is worked out.

\begin{sledst} \label{importpositive} We have $c_A^\Gamma\geqslant 0$. Moreover, $\Gamma$ is important if and only if $c_A^\Gamma>0$.
\end{sledst} To deduce positivity from importance here, apply Lemma \ref{lequiv_rel} to $\Psi_i$ and $\Psi$. 
Deducing importance from positivity is straightforward (we do not need this implication in what follows).
\begin{sledst} \label{sethyp} If $\dim E=-k<0$ for an essential facing $E$ of a relevant tuple $A$, then there exist essential facings $E_i$ of dimension $-i$ for $i=1,\ldots,k$, such that we have $E=E_k$, and $E_{i}$ is adjacent to $E_{i-1}$ for every $i$.
\end{sledst} To prove it, apply Theorem \ref{thess_main} by induction on $\dim E$.
Recall that adjacency is not transitive, and we cannot claim that $E_k$ is adjacent to $E_1$.

We now apply these results to study resultant sets.  Recall that we denote $\bigoplus_i \C[A(i)]$ by $\C[A]$, and the \textit{resultant set}
$R_A$ is the closure of the set of all $f\in\C[A]$ such that the system of equations $f=0$ has solutions in $\CC^n$.

\noindent\hspace{1cm}{\sc Notation.}  If $\Gamma$ is a face of $AI$, and $f\in\C[A]$ is
a tuple $\sum_{a\in A(i)}c_{a,i}x^a,\, i=0,\ldots,k,$ then the tuple $\sum_{a\in \Gamma(i)}c_{a,i}x^a,\, i\in I,$ in $\C[\Gamma]$ is denoted by $f^\Gamma$.
For every set $S\subset\C[\Gamma]$, we denote the set $\{f\, |\, f^\Gamma\in S\}\subset\C[A]$ by the same symbol $S$ (for example, $R_\Gamma\subset\C[A]$ is the set of all $f\in\C[A]$ such that the system $f^\Gamma=0$ has a solution).
\begin{lemma} \label{lunivres} 1) We have $\codim R_A=-\min\dim A$.
\newline 2) The resultant set $R_A$ is irreducible.
\newline 3) If $B$ is a face of $A$, then $R_B\subset R_A$.
\newline 4) For a generic $f\in R_A$, we have $\dim\{f=0\}=-\min\dim A-1$, and the Euler characteristic of $\{f=0\}$ is non-zero, unless $\sum A$ is contained in an affine hyperplane.
\end{lemma}
\textsc{Proof.} Part 1 is proved in \cite{sturmf}, Parts 2 and 3 are well known.

2 and 3. Let $T$ be a toric variety, whose fan is compatible with the convex hull of $\sum A$. Consider elements of $\C[A(i)]$ as sections of the linear bundle on $T$, corresponding to the convex hull of $A$ (see e.g. \cite{kh77}). The set $\{(f,x)\, |\, f(x)=0\}\subset\C[A]\times T$ is irreducible, because it is the total space of a vector bundle over $T$. Its projection to $\C[A]$ equals the union of $R_B$ over all faces $B\prec A$. Since $\bigcup_{B\prec A}R_B$ is irreducible, and $\codim R_A=-\min\dim A\leqslant-\min\dim B=\codim R_B$ for every $B$, then $R_A=\bigcup_{B\prec A}R_B$.

4. The equality $\dim\{f=0\}=-\min\dim A-1$ for generic $f$ follows from the equalities $\dim\{(f,x)\, |\, f(x)=0\}=|A|-1$ and $\dim R_A=|A|+\min\dim A$. With no loss in generality, assume that $A(0),\ldots,A(p)$ is the maximal essential subtuple of $A$. Denote the vector space, generated by pairwise differences of points of $A(0)+\ldots+A(p)$, by $L$. For every $f\in\C[A]$, the set $\{f_0=\ldots=f_p=0\}$ can be represented as $T\cdot C_f$, where $T\subset\CC^n$ is the subtorus, whose tangent plane at $1\in\CC^n$ equals $L^\perp\subset (\Z^n)^*$, and $C_f\subset\CC^n/T$ is an algebraic set. Applying the first statement of Part 4 to $R_{A(0),\ldots,A(p)}$, we conclude that, for generic $f\in R_A$, the set $C_f$ is non-empty and finite. Moreover, for every $c\in C_f$ and generic $f_i,\, i>p$, the restrictions
of $f_i,\, i>p$, to $T\cdot c$ define a nondegenerate complete intersection $V_c$. By the Kouchnirenko-Bernstein-Khovanskii formula \cite{kh77}, denoting the convex hull of the image of $A(i)$ in $\R^n/L$ by $B(i)$, the Euler characteristic of $V_c$ equals
$$\chi=(-1)^{\codim L+p-k} \sum\limits_{a_{p+1}+\ldots+a_k=\codim L}B(p+1)^{a_{p+1}}\cdot\ldots\cdot B(k)^{a_{k}}.\eqno{(*)}$$
Since $A_0,\ldots,A_p$ is the maximal essential subtuple of $A$, then
the polytopes $B(p+1),\ldots,B(k)$ in $\R^n/L$ are linearly independent, then, by Lemma \ref{lequiv1},
at least one of the terms on the right hand side of $(*)$ is non-zero, then $e\{f=0\}=|C_f|\cdot\chi$ is also non-zero of sign $(-1)^{\codim L+p-k}$. $\quad\Box$

\begin{sledst}\label{reshyp}
1) If $AI$ is the maximal essential subtuple of $A$, then $R_A=R_{AI}$.
\newline 2) If $\min\dim B<\min\dim A=0$ 
for a face $B\prec A$, then there exists another face $C\prec A$, such that $R_C\supset R_B$ is a hypersurface.
\end{sledst}
\textsc{Proof.}
1. We have $\codim p^{-1}(R_{AI})=\codim R_A=\min\dim A$ by Lemma \ref{lunivres}.1, and $p^{-1}(R_{AI})\supset R_A$ by definition.
Since both are irreducible by Lemma \ref{lunivres}.2, the converse inclusion follows.

2. Applying Corollary \ref{sethyp} with $E$ equal to the essential facing of $B$, we have (in the notation of Corollary):
$R_B=R_E=R_{E_k}\subset\ldots\subset R_{E_1}$ by Lemma \ref{lunivres}.3, and $\codim R_{E_1}=-\min\dim E_1=1$. $\quad\Box$

\begin{rem} Let us say that an essentual facing $\Gamma$ of $A$ is \textit{essentially adjacent} to an essential facing $\Delta$,
and write $\Gamma\preccurlyeq\Delta$, if $R_\Gamma\subset R_\Delta$. Then we have:
\newline 1) $\preccurlyeq$ is transitive,
\newline 2) $\Gamma\preccurlyeq\Delta\; \Rightarrow\; \dim \Gamma\leqslant\dim\Delta$ (by Lemma \ref{lunivres}),
\newline 3) For every $\Gamma_1\preccurlyeq\Gamma\preccurlyeq\Gamma_2$, such that $\dim\Gamma_1<\dim\Gamma<\dim\Gamma_2$,
there exist essential facings $\Delta_1$ and $\Delta_2$, such that
$$\Gamma_1\preccurlyeq\Delta_1\preccurlyeq\Gamma\preccurlyeq\Delta_2\preccurlyeq\Gamma_2 \mbox{ and } \dim\Delta_1+1=\dim\Gamma=\dim\Delta_2-1$$ (by Theorem \ref{thess_main}).

In this sense, essential adjacency of essential facings of a tuple of polytopes behaves in the same way
as adjacency of faces of a polytope. However, essential adjacency is weaker than the adjacency of Definition \ref{dessf},
because the latter is not transitive. It would be interesting to give a combinatorial definition of
essential adjacency, not relying upon resultant sets (and thus applicable to non-integer polytopes as well), e.g. the following one.
\end{rem}
\begin{conj} \label{conjess_main} Essential adjacency is the transitive closure of adjacency of Definition \ref{dessf}
(in particular, if $\Gamma\preccurlyeq\Delta$ and $\dim\Gamma=\dim\Delta$, then $\Gamma=\Delta$).
\end{conj}
If $\min\dim A=0$, then $R_A=\C[A]$ is trivial, and, instead of it, we study the set
$$D_A=\{f\, |\,  f(x)=0 \mbox{ and } \sum\lambda_i df_i(x)=0 \mbox{ for some }\lambda_i\ne 0 \mbox{ and } x\in\CC^n\}\subset\C[A].$$
This set will help us to study the desired set of degenerate systems
$$B_A=\{f\, |\, 0 \mbox{ is a critical value of } f^\Gamma:\CC^n\to\C^{|\Gamma|} \mbox{ for some } \Gamma\}\subset\C[A],$$
and is simpler than $B_A$ (for example, it is irreducible, while $B_A$ is not, and all irreducible components of $B_A$ equal $D_\Gamma$ for certain $\Gamma$, see below for a short account and \cite{dcg} for details). See also \cite{st11} for more information on the set $D_A$ for $k=n-1$, where its closure is called the \textit{mixed discriminant variety}. The following is well known.


\begin{lemma} \label{ld111} 1) We have $f\in D_A \Leftrightarrow f_{\{0,\ldots,k\}}\in D_{A_{\{0,\ldots,k\}}}$ (in the notation, introduced for Theorem \ref{lozungdecomp}).
\newline 2) The set $D_A$ is irreducible.
\newline 3) If $k=0$, pairwise differences of points of $A$ generate $\Z^n$, and $D_A$ is a hypersurface, then, for a generic polynomial $f\in D_A$, the critical point of the map $f:\CC^n\to \C^{k+1}$ in $f^{-1}(0)$ is unique and nondegenerate.
\end{lemma}
{\sc Proof.} Part 1 is a tautology, and reduces Part 2 to the case $k=0$. In this case, consider the set $S=\{ (x,f) \, |\, f(x)=df(x)=0\}\subset\CC^n\times\C[A]$. The projection $S\to\CC^n$ is a vector bundle of rank $|A|-\dim(\mbox{convex hull of } A)-1$, thus $S$ is irreducible of dimension $|A|+n-\dim(\mbox{convex hull of } A)-1$,
thus its projection $D_A\subset\C[A]$ is irreducible, which proves Part 2. 

Moreover, if the convex hull of $A$ is $n$-dimensional, and $D_A$ is a hyprsurface, then $\dim D_A=\dim S$, and the projection $S\to D_A$ is generically finite. By the Sard lemma, a generic fiber of this projection is finite and nondegenerate. If $(x,f)$ is a point in such a fiber, then a straightforward computation shows that $f$ is a smooth point of $D_A$, and the covector $\sum_{a\in A} x^a dc_a$ is normal to $T_f D_A$  (recall that $(c_a,\, a\in A)$ is the natural coordinate system on $\C[A]$ such that $c_a(f)$ is the coefficient of the monomial $x^a$ in $f\in\C[A]$). This condition uniquely defines $x$ for a given $f$ in the smooth part of $D_A$, thus the projection $S\to D_A$ is a bijection over the smooth part of $D_A$. $\quad\Box$

\begin{theor}[\cite{dcg}, Corollary 2.32] \label{thdcg} 
The codimension 1 part of $B_A$ equals the union of codimension 1 sets $D_\Gamma$ over all important $\Gamma$ (Definition \ref{defimport}).
\end{theor}
Let $i_\Gamma$ be the index of the lattice, generated by the pairwise differences of the points of $\sum\Gamma$
(recall that the index of $L\subset\Z^n$ is the cardinality of the set $(L\otimes\Q\cap\Z^n)/L$). Recall that we denote the Euler characteristic by $e$.
\begin{theor} \label{theuler} If $D_\Gamma\subset\C[A]$ is a hypersurface, then, for generic $\tilde f$ in it and generic $f$ in $\C[A]$,
we have $e\{\tilde f=0\}-e\{f=0\}=(-1)^{n-k}i_\Gamma c^\Gamma_A$.
\end{theor}
The proof is given below. From these two Theorems and Corollary \ref{importpositive} we get
\begin{sledst} \label{coinc1} For generic systems $f$ in $\C[A]$ and $\tilde f$ in the codimension 1 part of $B_A$,
the difference $e\{\tilde f=0\}-e\{f=0\}$ equals $(-1)^{n-k}i_\Gamma c^\Gamma_A$, and, in particular, is non-zero of sign $(-1)^{n-k}$.
\end{sledst}
In the notation, introduced for Theorem \ref{lozungdecomp}, the set $\{f_I=0\}\subset\CC^I\times\CC^n$ is invariant
under the action of $\CC$ on $\C^I$ by dilatations, and we denote the quotient by $\{f_I=0\}/\CC$.
\begin{lemma} \label{l241}
1) For every tuple $f$ of Laurent polynomials on $\CC^n$, we have $e\{f=0\}=\sum_I e\{f_I=0\}/\CC$.
\newline 2) For every face $\Gamma$ of $AI$, we have $c^{\Gamma_I}_{A_{\{0,\ldots,k\}}}=\sum_{J\supset I}c^{\Gamma}_{AJ}$.
\end{lemma}
See \cite{dcg}, Lemma 1.7, for Part 2, and the equality of the Euler characteristics in its proof for Part 1.

\textsc{Proof of Theorem \ref{theuler}.} By Lemma \ref{l241}, Lemma \ref{ld111}(1) and the equality $i_\Gamma=i_{\Gamma_{\{0,\ldots,k\}}}$ we can reduce the statement to the case $k=0$,
studying the hypersurfaces $\{f_I=0\}/\CC$ instead of the complete intersection $\{f=0\}$.
Assuming $k=0$, let $\widetilde M$ and $M$ be the closures of the hypersurfaces $\{\tilde f=0\}$ and $\{f=0\}$
in the $A$-toric variety $T\supset\CC^n$. Then $M$ is transversal to all orbits of $T$, and, by Lemma \ref{ld111}(3), so
is $\widetilde M$ except for one nondegenerate singular point $x$ at the $\Gamma$-orbit. If $f$ is
close to $\tilde f$, then $M$ is diffeomorphic to $\widetilde M$ outside a small neighborhood $U$ of $x$,
and $M\cap U$ is the Milnor fiber of $\widetilde M$ at $x$. Thus, the desired difference $e\{\tilde f=0\}-e\{f=0\}$
equals $(-1)^{n-k}$ times the Milnor number of $\widetilde M$ at $x$. Passing to the normalization $T$,
we have $i_\Gamma$ points over $x$, 
and the Milnor number at each
of them equals $c^\Gamma_A$ (\cite{dcg}, Theorem 1.18). $\quad\Box$

\section{Discriminant}

In Section \ref{Seuldiscr}, we introduce the Euler discriminant for an arbitrary morphism of complex algebraic varieties (Definition \ref{defeuldiscr}). For proper maps of smooth varieties, the Euler discriminant turns out to be invariant under transversal base changes (Theorem \ref{thbasech}), which is actually a matter of differential topology and can be extended to real manifolds.

In Section \ref{Sbertdiscr}, we define the Bertini discriminant for a fiberwise compactification of a morphism of complex algebraic varieties (Definition \ref{defbertdiscr}). If the domain of the compactified map is a local complete intersection, then the Bertini discriminant turns out to be a hypersurface under certain mild assumptions on the dimension of adjacent strata of the compactification (Theorem \ref{thgeneral}.4). The Bertini and the Euler discriminants estimate the bifurcation set of a morphism from above and from below respectively, and both of them turn out to be invariant under base changes, satisfying certain transversality assumptions (Theorem \ref{thgeneral}.1--3).

\subsection{Euler discriminant} \label{Seuldiscr} Let $\pi:Q\to P$ be a morphism of complex algebraic varieties. Whenever we make a statement about a \textit{germ} of $\pi$ near $p\in P$, we imply the existence of a complex neighborhood of $p$, such that the statement is valid for the restriction of $\pi$ to the preimage of every smaller neighborhood. The \textit{bifurcation set} of $\pi$ is the minimal subset $B_\pi\subset P$ such that the restriction of $\pi$ to the preimage of $P\setminus B_\pi$ is a locally trivial fibration. 
By the Hardt theorem, one can decompose $P$ into smooth strata $P_i$, such that the Euler characteristic of the fiber $\pi^{-1}(x)$ equals the same number $e_i$ for all $x\in P_i$. Let $P_0$ be the dense stratum, and $P_1,\ldots,P_k$ be the codimension 1 strata. \begin{defin}\label{defeuldiscr} The Weil divisor $(-1)^{\dim Q-\dim P}\sum_{i=1}^k (e_i-e_0)\overline{P}_i$ is called the \textit{Euler discriminant} of $\pi$ and is denoted by $E_\pi$. \end{defin} 
For morphisms $\pi:Q\to P$ and $f:R\to P$ of complex algebraic manifolds, the \textit{fiber product} $Q\times_PR$ is the intersection of (graph of $f)\times Q$ and $R\times($graph of $\pi)$ in $R\times P\times Q$, and $f$ is said to be \textit{transversal} to $\pi$, if the smooth submanifolds (graph of $f)\times Q$ and $R\times($graph of $\pi)$ are transversal.
We recall some standard facts about transversality:
\begin{utver} \label{proptransv} 1. If $f:R\to P$ and $\pi:Q\to P$ are transversal, then $Q\times_PR$ is also a manifold, and $f(R)$ is not contained in the set of critical values of $\pi$. \par 2 (transitivity). If, in addition, $g:S\to R$ is transversal to the projection $Q\times_PR\to R$, then the composition $f\circ g$ is transversal to $\pi$. 
\end{utver}
\textsc{Proof.} The fiber product $Q\times_PR$ is smooth by definition. Without loss of generality, transitivity can be proved for linear maps $f, g, \pi$. This special case is obvious: if $\im\pi+\im f=P$ and $\im g+f^{-1}(\im\pi)=R$, then $\im fg+\im\pi=P$. By transitivity, the fact that $f(R)$ is not contained in the discriminant of $\pi$ can be reduced to the case of $\dim R=0$, which follows by definition. $\quad\Box$\par Recall that, for a smooth map $f:R\to P$ and a divisor $D$ in $P$, the pull-back $f^*D$ is a divisor in $R$ that additively depends on $D$ and is defined near $x\in R$ by the following property: if $D$ is locally irreducible near $f(x)$ and given by an equation $g=0$, then $f^*D=0$ whenever $g\circ f=0$ identically, and $f^*D$ is given by the equation $g\circ f=0$ otherwise.
\begin{theor}[base change]\label{thbasech} If $f:R\to P$ and a proper $\pi:Q\to P$ are transversal, then the Euler discriminant of the projection $Q\times_PR\to R$ equals $f^* E_\pi$.
\end{theor}
\textsc{I. Proof for a germ $f:(\C,0)\to(P,p)$.} If $f$ is a germ of a smooth curve, transversal to the Hardt strata $P_i$ of $P$, then the statement is valid by definition of the Euler discriminant. In general, the image of $f$ intersects the discriminant of $\pi$ at an isolated point $p$ by Proposition \ref{proptransv}(1), and we can consider its perturbation $F:(\C^2,0)\to(P,p),\, F(x,0)=f(x)$, such that $F$ is transversal to $\pi$, and $F(\cdot,t)$ is transversal to the Hardt strata $P_i$ for small $t\ne 0$. More precisely, we choose the product $D\times T\subset\C^2$ of two open disks, such that the restriction $f_t$ of the map $F$ to the disk $D_t=D\times\{t\}$ is transversal to every $P_i$ for every non-zero $t\in T$, both $Q\times_P(D\times T)$ and $Q\times_PD_0$ retract to the fiber $\pi^{-1}(p)$, and $E_\pi\cap (D\times T)$ retracts to $p$. By transversality of $F$ and $\pi$, the homeomorphism type of $Q\times_PD_t$ does not depend on $t$, thus $$e(Q\times_PD_0)=e(Q\times_PD_t)\mbox { for } t\in T.$$ Denoting the total multiplicity of a 0-dimensional divisor by $m$, we have $$m f_0^* E_\pi=m f_t^* E_\pi \mbox { for } t\in T.$$ By transversality of $f_t$ to the Hardt strata, $f_t(D)$ intersects certain hypersurfaces $P_i,\, i\in I$, transversally at isolated points $p_i$, and does not intersect other strata, thus we have $$(-1)^{\dim Q-\dim P}m f_t^*E_\pi=\sum_{i\in I}(e\pi^{-1}(p_i)-e_0),$$ where $e_0$ is the Euler characteristic of a generic fiber of $\pi$. By additivity of the Euler characteristic, the latter sum equals $$\sum_{i\in I}(e\pi^{-1}(p_i)-e_0)=e(Q\times_PD_t)-e_0.$$ As a result, we have $$(-1)^{\dim Q-\dim P}m f_0^*E_\pi=(-1)^{\dim Q-\dim P}m f_t^*E_\pi=\sum_{i\in I}(e\pi^{-1}(p_i)-e_0)=$$ $$=e(Q\times_PD_t)-e_0=e(Q\times_PD_0)-e_0=e(\pi^{-1}(p))-e_0,$$ which implies the desired statement.\par
\textsc{II. Proof for arbitrary $f$.} Part I allows us to reformulate the definition of the Euler discriminant as follows. The Euler discriminant of an arbitrary proper map $\pi:Q\to P$ is the (unique) Weil divisor $D$ in $P$, such that, for every germ $j:(\C,0)\to(P,p)$, transversal to $\pi$, the intersection number $m j^*D$ is equal to $(-1)^{\dim Q-\dim P}(e\pi^{-1}(p)-e_0)$, where $e_0$ is the Euler characteristic of a generic fiber of $\pi$. This definition is invariant under transversal base changes by Proposition \ref{proptransv}(2). $\quad\Box$\par
The study of Euler discriminants of complete intersections can be reduced to the case of hypersurfaces by the following version of the Cayley trick.
\begin{theor} \label{theulcayley} Let $\pi$ be the restriction of a projection  $M\to P$ to a complete intersection $\{f_0=\ldots=f_k=0\}\subset M$. Let $\lambda_0:\ldots:\lambda_k$ be the standard coordinates on $\CP^k$, and let $\pi_0$ be the restriction of the projection $\CP^k\times M\to M\xrightarrow{\pi}P$ to the set $\{\sum_i\lambda_if_i=0\}\subset\CP^k\times M$. Then we have $E_{\pi_0}=E_\pi$.
\end{theor}
\textsc{Proof.} By additivity of the Euler characteristic, for every $y\in P$ we have $e\pi^{-1}(y)=e\pi_0^{-1}(y)-ke(M)$. The constant term $ke(M)$ does not affect the definition of the Euler discriminant.
$\quad\Box$\par

\subsection{Bertini discriminant} \label{Sbertdiscr} 
Let $\pi$ be a morphism of a complex quasi-projective scheme $Q$ to a smooth variety $P$. The following notions are well known, and their various versions were studied by many authors mentioned in the introduction.
\begin{defin} \label{defbertdiscr} A proper morphism $\bar\pi:\bar Q\to P$ of a complex quasi-projective scheme $\bar Q$, endowed with a collection of Cartier divisors $D_i$, is called a \textit{properization} of $\pi$, if $Q$ is an open subscheme of $\bar Q$ with $\pi=\bar\pi|_Q$, and, set-theoretically, $Q=\bar Q\setminus \bigcup_i D_i$. A point $y\in P$ is said to be in the \textit{singular locus} $\sing\bar\pi$, if it equals $\bar\pi (x)$ such that
\newline 1) $x\in\sing \bar Q$, or
\newline 2) the differentials of the local equations of the divisors $D_i$, passing through $x$, are linearly dependent at $x$. \newline
A point $y\in P$ is said to be in the \textit{Bertini discriminant} $\bar B_{\bar \pi}$ , if it equals $\bar\pi (x)$ such that
\newline 1) $x\in\sing \bar Q$, or
\newline 2) the differential $d\bar\pi$ is not surjective at $x$, or
\newline 3) the differentials of the local equations of the divisors $D_i$, passing through $x$, have linearly dependent restrictions to $\ker d\bar\pi$ at $x$. \newline
The properization $\bar\pi$ is said to be \textit{c-smooth}, if $\codim\sing\bar\pi>c$. For a map $F:R\to P$, the \textit{fibered product} $R\times_P \bar Q$ is the (scheme-theoretic) intersection of $R\times($graph of $\bar\pi)$ and $($graph of $F)\times \bar Q$ in $R\times P\times \bar Q$. The \textit{induced properization} $F^*\bar\pi$ is the projection of $R\times_P \bar Q$ to $R$, endowed with the collection of Cartier divisors $R\times_P D_i$. \end{defin} 1-smooth properizations of $\pi$ provide convenient upper and lower bounds for the bifurcation set $B_\pi$ in the sense of the following theorem.
For every set $I$, denote the (scheme-theoretical) intersection $\bigcap_{i\in I}D_i$ by $\bar Q_I$ and its open subscheme $\bigcap_{i\in I}D_i\setminus\bigcup_{j\notin I}D_j$ by $Q_I$. In particular, $Q_\varnothing=Q$ and $\bar Q_\varnothing=\bar Q$.

\begin{theor} \label{thgeneral} 1. We have $|E_\pi|\subset B_\pi\subset\bar B_{\bar\pi}$.
\newline
2. If $\bar\pi$ is 0-smooth, then $\bar B_{\bar\pi}$ is nowhere dense.
If $\bar\pi=F^*\bar p$, then we have $\bar B_{\bar\pi}=F^{-1}\bar B_{\bar p}$.
\newline
3. If, moreover, both $\bar\pi$ and $\bar p$ are 1-smooth, then we have $E_{\pi}=F^*E_p$.
\newline
4. If $\bar\pi$ is 0-smooth and $\bar Q$ is a local complete intersection, then we have $\codim\bar B_{\bar\pi}=1$, provided that

a) For every $I$ such that $\bar\pi(\bar Q_I)$ is of positive codimension, there exists $J$ such that $\bar\pi(\bar Q_J)$ has
codimension 1 and contains $\bar\pi(\bar Q_I)$.

b) For every $I$, the set $Q_I$ is affine.\end{theor} 
We will use it in the context where $Q$ is a complete intersection in a complex torus, $\pi$ is its projection to a subtorus, and $\bar Q$ is its smooth fiberwise toric compactification, in order to prove that generically both inclusions of Part (1) are equalities.

\textsc{Proof of Part 1.} The first inclusion means that all fibers of a locally trivial fibration have the same Euler characteristic.
The second inclusion follows from the fact that, outside the Bertini discriminant $\bar\pi^{-1}(\bar B_{\bar\pi})$,
the projection $\bar\pi$ defines a stratified locally trivial fibration of the stratified space $\bar Q=\bigsqcup_I Q_I$.

\textsc{Proof of Part 2.} 
The difference $\bar B_{\bar\pi}\setminus\sing\bar\pi$ is nowhere dense by the Bertini-Sard theorem, and so is $\sing\bar\pi$ by 0-smoothness.
Let $\hat F$ be the induced map of the domain of $\bar \pi=F^*\bar p$ to the domain of $\bar p$. Since $\hat F$ induces the isomorphism between the
(scheme-theoretic) fibers $\bar\pi^{-1}(y)$ and $\bar p^{-1}(F(y))$, it preserves Conditions (1) and (2) in the definition of the Bertini discriminant. Moreover, if these fibers are regular at some point $x$, then $\hat F$ preserves the restrictions of the local equations of the divisors $D_i$, restricted to $T_x\bar\pi^{-1}(y)$, and, in particular, preserves Condition (3) in the definition of the Bertini discriminant.  

\textsc{Proof of Part 3.} Denoting the restriction of $\bar\pi$ to $\bar Q_I$ by $\pi_I$ (so that $\pi_\varnothing=\bar\pi$), we have $E_\pi=\sum_I(-1)^I E_{\pi_I}$
by the inclusion-exclusion formula and additivity of Euler characteristic. For the restriction of $E_{\pi_I}$ to $P\setminus\sing\bar\pi$,
we have the desired equality $E_{\pi_I}=F^*E_{p_I}$ for every $I$ by Theorem \ref{thbasech}, because $F$ and $p_I$ are transversal outside $\sing\bar\pi$. Since $\codim\sing\bar\pi>1$, we have $E_{\pi_I}=F^*E_{p_I}$ on $P$ as well, and apply it to every term of the sum $E_\pi=\sum_I(-1)^I E_{\pi_I}$.

\textsc{Proof of Part 4.}
We prove that $\overline{B}_{\bar\pi}$ has a codimension 1 component at its arbitrary point $y$.

\textbf{I.} Let $\Sigma_{I,y}$ consist of all points $x\in\bar\pi^{-1}_I(y)$, such that $\bar Q$ is not regular at $x$, or the differential $d\bar\pi$ is not surjective at $x$, or the differentials of the local equations of $D_i,\, i\in I,$ are linearly dependent on $\ker d\bar\pi$.
Choose a maximal $I$ such that $\Sigma_{I,y}$
is not empty. Since $\Sigma_{I,y}\cap D_j\subset \Sigma_{I\cup\{j\},y}$ for every $j$, the maximality of $I$ implies that $\Sigma_{I,y}\cap D_j=\varnothing$ for all $j\notin I$, i. e. the compact set $\Sigma_{I,y}$ is contained in $Q_I$. Since the latter set is affine (by Condition b), the set $\Sigma_{I,y}$ is finite.

\textbf{II.} Without loss of generality, we can assume that $|I|\leqslant\dim Q-\dim P$. Otherwise
$y$ is contained in the set $\bar\pi(\bar Q_I)$ of positive codimension, thus (by Condition a) it is also contained in $\bar\pi(\bar Q_J)$ of codimension $1$, which is in turn contained in $\bar B_{\bar\pi}$.

\textbf{III.} Near a point $x\in \Sigma_{I,y}$, the local complete intersection $Q_I$ is given by local equations $h_1=\ldots=h_k=0$ in an ambient space $C^N$, the divisor $D_i$ is given by the equation $g_i=0$,
and the fiber $\bar\pi^{-1}_I(y)$ is given by the equations $u_1\circ\bar\pi=\ldots=u_p\circ\bar\pi=0$, where $(u_1,\ldots,u_p)$ is a system of local coordinates near $y\in P$. The system of equations $$h_j=0,\, j=1,\ldots,k,\, g_i=0,\, i\in I,\,
u_l\circ\pi=\varepsilon_l,\, l=1,\ldots,p,\eqno{(*)}$$ defines
a family of schemes $C_\epsilon$, parameterized by $\epsilon=(\varepsilon_1,\ldots,\varepsilon_p)$. The scheme $C_0$ has an isolated complete intersection
singularity, because $x$ is isolated in $\Sigma_{I,y}$ by Step (I), and the number of equations in $(*)$ is not greater than $N$ by Step (II).
Thus the family $C_\epsilon$ is induced from the versal deformation of the isolated complete intersection singularity $C_0$, then its discriminant $\{\epsilon\, |\, C_\epsilon \mbox{ is singular } \}\subset \bar B_{\bar\pi}$, which contains $y$, is a hypersurface, because the discriminant of the versal deformation of icis is a hypersurface. $\quad\Box$

\begin{rem} If every fiber of $\pi$ contains finitely many critical points and has the same dimension, then the discriminant of $\pi$ in the sense of Teissier is defined (see \cite{teissier}), and the Euler discriminant $E_\pi$ coincides with its codimension 1 part. 
It would be interesting to construct a common generalization of the Euler discriminant and the codimension 1 part of the Teissier discriminant, so that it would be defined under relaxed assumptions on both 1-smoothness of $\pi$ and finiteness of its restriction to the set of critical points, and still invariant under base changes.\par  For instance, if we admit $Q$ with isolated complete intersection singularities and $P$ only 1-dimensional, then, in Definition \ref{defeuldiscr}, we should replace the Euler characteristic of the fiber $\pi^{-1}(x)$ with the sum of the Euler characteristics of the set $\pi^{-1}(x)\setminus\sing Q$ and of the local Milnor fibers of $Q$ at the points $\pi^{-1}(x)\cap\sing Q$. This version of the Euler discriminant is invariant under arbitrary 1-dimensional base changes.

In the same way, using integration with respect to the Euler characteristic, one can define $E_\pi$ for $Q$ with non-isolated singularities. E.g. assume that $Q$ is a hypersurface with arbitrary singularities, and denote the Euler characteristic of the local Milnor fiber of $Q$ at $x$ by $\mu(x)$. Define $E_\pi$, replacing in its Definition \ref{defeuldiscr} the Euler characteristic of the fiber $\pi^{-1}(x)$ with the Euler characteristic integral of the function $\mu:Q\to\Z$ over the fiber. Then $E_\pi$ is still invariant under 1-dimensional base changes, however this generality is already irrelevant to the topic of our paper: the example constructed in \cite{nemethicanadian} shows that we cannot anymore expect that, generically, the inclusions of Theorem \ref{thgeneral}.1 are equalities. \end{rem}

\section{Toric case}
We apply combinatorics of Section 2 and topology of Section 3 to study the Euler and the Bertini discriminants of a projection of a nondegenerate complete intersection in a complex torus, and find that the two versions of the discriminant coincide in this case. This yields all the facts, announced in the introduction: see Section \ref{Sunivcase} for Theorems \ref{thnondeg}--\ref{lozungdecomp}, Section \ref{Sgencaseirr} for Proposition \ref{lozungirr} and Section \ref{Sgencaserel} for Theorem \ref{lozung} and Proposition \ref{lozungnewt}.

\subsection{Universal case}\label{Sunivcase} We now prove Theorems \ref{thnondeg}--\ref{lozungdecomp}. All of them are corollaries of the following one.
Let $Q\subset \CC^n\times\C[A]$ consist of all pairs $(x,f)$ such that $f(x)=0$, and denote the projection of $Q$ to $\C[A]$ by $p$.
Choose a smooth toric compactification $T$ of $\CC^n$, such that its fan $\Gamma$ is compatible with the convex hulls of $A_0,\ldots,A_k$,
denote the closure of $Q$ in $T\times\C[A]$ by $\bar Q$, the projection of $\bar Q$ to $\C[A]$ by $\bar p$, and the intersections of $\bar Q$ with the closures of the codimension 1 orbits of $T$ by $D_1,D_2,\ldots$. Consider the Euler discriminant $E_p$ (Definition \ref{defeuldiscr}) and the bifurcation set $B_p$ of the projection $p$, and the Bertini discriminant $\bar B_{\bar p}$ (Definition \ref{defbertdiscr}) of the properization $(\pi,D_1,D_2,\ldots)$. 
\begin{theor}\label{thmainuniv} If $A$ is relevant, then $E_p=E_A$ (Definition \ref{defaeuldiscr}), and $|E_p|=B_p=\bar B_{\bar p}$ (in particular, all of them are hypersurfaces).
\end{theor}
\textsc{Proof.} The equality $E_p=E_A$ follows by definition. By Corollary \ref{reshyp}.2, the properization $(p,D_1,D_2,\ldots)$ satisfies Condition 3a of Theorem \ref{thgeneral}, the other conditions
are satisfied by definition, thus, by Theorem \ref{thgeneral}.4, $\bar B_{\bar p}$ is a hypersurface. By Corollary \ref{coinc1}, every codimension 1
component of $B_A=\bar B_{\bar p}$ participates in $E_p$ with a non-zero multiplicity, thus $|E_p|=\bar B_{\bar p}$. The remaining equality follows by Theorem \ref{thgeneral}.1. $\quad\Box$.

We now use this to deduce the three theorems of the introduction.

\textsc{Proof of Theorem \ref{thnondeg}.} If $\min\dim A<0$, then the condition $f\in R_A$ is obviously equivalent to (1) and (3), and is equivalent to (2) by Lemma \ref{lunivres}.4. Otherwise, by definition of the Bertini discrimiant,
$\bar B_{\bar p}$ consists of all $f$ that satify Condition (3).

If $f$ is in $\bar B_{\bar p}$, then it satisfies (2), because we have arbitrarily small $\tilde f$ and $\tilde f'$, such that $f+\tilde f$ is a generic point of $|E_p|$ and $f+\tilde f'$ is outside $|E_p|$, thus the sets $\{f+\tilde f=0\}$ and $\{f+\tilde f'=0\}$ have different Euler characteristic by definition of $E_p$, thus at least one of them differs from $e\{f=0\}$. Since $f$ satisfies (2), it also satisfies (1).

If $f$ is not in $\bar B_{\bar p}$, then it does not satisfy (1) and (2), because it is not in $B_p$. $\quad\Box$

\textsc{Proof of Theorem \ref{lozung0}.} $B_p$ is the support set of the divisor $E_p$. $\quad\Box$

\textsc{Proof of Theorem \ref{lozungdecomp}.} If $k=0$, then, by Theorem \ref{theuler}, the left hand side defines the divisor
$\sum_{\Gamma} i_\Gamma c^\Gamma_A D_\Gamma$ with $\Gamma$ running over all faces of $A$, and the right hand side defines
the same divisor by the prime decomposition theorem \cite{gkz}, Chapter 10, Theorem 1.2. The case of arbitrary $k$ reduces to $k=0$
by means of the equality $E_A=\prod_{I\subset\{0,\ldots,k\}}E_{A_I}^{(-1)^{k+1-|I|}}$, which follows from Theorem \ref{theulcayley}. $\quad\Box$

\subsection{General case for irrelevant $A$}\label{Sgencaseirr} We come back to the setting and notation of Section \ref{Sintrobif}. We restrict our attention to the setting (I. Laurent polynomials), the other one can be treated in the same way. Recall that $\min\dim A$ is the minimum of $p-\dim(\mbox{convex hull of } A_{i_1}+\ldots+A_{i_p})$ over all sequences $i_1<\ldots<i_p$ (including the empty one).

\begin{utver} \label{irrelpure} 1) If $\min\dim A\leqslant-c<0$, and $f$ is $(c-1)$-nondegenerate (Definition \ref{cdeg}), then the generic fiber of the projection $\pi:Q\to\CC^m$ is empty, the bifurcation set is the closure of $\pi(Q)$, and its codimension is at least $c$. \newline 2) If $\min\dim A=-c<0$, and $f$ is $c$-nondegenerate, then, moreover, the codimension is equal to $c$.
In this case, if the sum $A_0+\ldots+A_k$ is not contined in an affine hyperplane, then a generic non-empty fiber has a non-zero Euler characteristic.
\end{utver} 
If we empose genericity assumptions on $f$ as generously as we wish, then the statement
of Part 2 is almost obvious. The aim of this proposition is to provide a strict genericity assumption. Although one can deduce it from the results of the previous sections (or their versions for $\min\dim A<-1$) in the same way as it is done for relevant $A$ below, we prefer to give a direct (almost) self-contained proof, illustrating some ideas behind the complicated case of relevant $A$.
\begin{lemma} \label{tropcover} Let $\Sigma\subset\CC^n\times\CC^n$ be an algebraic variety, such that $\codim \pi(\Sigma)=c$. Then the tropical fan of $\Sigma$ contains
a linear function $v:\Z^n\oplus\Z^m\to\Z$ of dimension at most $c$ (Definition \ref{cdeg}).
\end{lemma}
\textsc{Proof.} Let $\int H\subset\Z^m$ be the fiber polytope (Definition \ref{defmp}) of the convex hull of $H_0+\ldots+H_m\subset\Z^n\oplus\Z^m$. Without loss of generality, we may assume that $\int H$ is not contained in a hyperplane. 
The tropical fan of $\pi(\Sigma)$ contains a liner function $u:\Z^m\to\Z$, whose support face in $\int H$ is at most $c$-dimensional, because the set of linear functions with higher-dimensional support faces has codimension $c-1<\codim \pi(\Sigma)$. Since the tropical fan of $\pi(\Sigma)$ is the projection of the tropical fan of $\Sigma$ (see \cite{kaz}, \cite{st}), $u$ can be extended to a linear function $v:\Z^n\oplus\Z^m\to\Z$, contained in the tropical fan of $\Sigma$. The dimension of $v$ is at most the dimension of $(\int H)^u$, which is at most $c$. $\quad\Box$

\textsc{Proof of Proposition.} \textbf{I.} If $\min\dim A\leqslant-c<0$, then $(c-1)$-nondegeneracy means that $\{f_0^v=\ldots=f_k^v=0\}$ is empty for every linear function $v:\Z^n\oplus\Z^m\to\Z$ of dimension smaller than $c$, i.e. the tropical fan of $Q$ does not contain linear functions of dimension smaller than $c$. Applying Lemma \ref{tropcover} with $\Sigma=Q$, we deduce $\codim \pi(Q)\geqslant c>0$. Thus, the generic fiber of $\pi:Q\to\CC^m$ is empty, and the closure of $\pi(Q)$ is the bifurcation set.

\textbf{II.} If, moreover, $f$ is $c$-nondegenerate, then, for generic affine linear functions $l_1,\ldots,l_{m-c}$ on $\CC^m$ and $l_{m-c+1},\ldots,l_{n+m-k-1}$ on $\CC^n\times\CC^m$, the system $$f_0=\ldots=f_k=l_1=\ldots=l_{n+m-k-1}=0 \eqno{(*)}$$ of $m+n$ equations of $m+n$ variables satisfies the genericity assumption of the Kouchnirenko-Bernstein theorem (\cite{bernst}). Thus the number of its solutions equals the mixed volume of the convex hulls of $H_0,\ldots,H_k$, $m-c$ copies of the standard simplex in $\Z^m$ and $c+n-k-1$ copies of  the standard simplex in $\Z^n\oplus\Z^m$. The condition $\min\dim A=-c$ implies that these polytopes are linearly independent, thus their mixed volume $V$ is positive (Lemma \ref{lequiv1}), thus the system $f_0=\ldots=f_k=l_1=\ldots=l_{m-c}=0$ is consistent for generic $l_1,\ldots,l_{m-c}$ on $\CC^m$. In other words, $\pi(Q)$ intersects a generic $c$-dimensional affine subspace in $\CC^m$, which implies that it has a component of codimension at most $c$. It remains to prove that it has no components of greater codimension.

\textbf{III.} Moreover,
denoting the plane $l_1=\ldots=l_{m-c}=0$ by $L\subset\CC^m$, 
consider the set $Z$ of pairs $(z,\tilde f)$, such that $\tilde f\in\C[H_0]\oplus\ldots\oplus\C[H_k]$, and $z$ is
in the intersection of $\pi\{x\, |\, \tilde f(x)=0\}$ and $L$. Then the projection of $Z$ to $\C[H_0]\oplus\ldots\oplus\C[H_k]$ is proper at $f$.
To prove it, consider the set $\hat Z$ of all pairs $(y,\tilde f)$, such that $\tilde f\in\C[H_0]\oplus\ldots\oplus\C[H_k]$, and
$y\in\CC^n\times\CC^m$ is a solution of the system $\tilde f_0=\ldots=\tilde f_k=l_1=\ldots=l_{n+m-k-1}=0$. The desired properness follows from the two observations below.

1) The projection of
$\hat Z$ to $\C[H_0]\oplus\ldots\oplus\C[H_k]$ is proper at $f$, because the system $(*)$ satisfies the genericity assumption of the Kouchnirenko-Bernstein theorem, which assures that the system $\tilde f_0=\ldots=\tilde f_k=l_1=\ldots=l_{n+m-k-1}=0$ has the same number of solutions for every $\tilde f$ near $f$.

2) The projection of $\hat Z$ to $Z$ is an epimorphism near $f$. To prove it at a point $(z,f)$, note that
the point $z$ is contained in a codimension $c$ component of $\pi(Q)$, because it has no components of lower codimension,
and components of higher codimension do not intersect the generic plane $L$.
Thus, the fiber of $Q\to \pi(Q)$ over $z$ is at least $(c+n-k-1)$-dimensional, because otherwise $\codim Q>k+1$ at this fiber,
which contradicts the definition $Q=\{f_0=\ldots=f_k=0\}$. Thus, the generic plane $l_{m-c+1}=\ldots=l_{n-m-k-1}=0$ of
complementary dimension intersects this fiber at some point $y$. Thus, $(z,f)\subset Z$ is the projection of $(y,f)\subset \hat Z$.

\textbf{IV.} Assume the contrary to what remains to prove: assume that the closure of $\pi(Q)$ has a component $R$ of codimension greater than $c$. Consider the set $S$ of all consistent systems of equations $h_0=\ldots=h_k=0,\, h_i\in\C[A_i]$, denote the closure of $S$ by $\bar S$, and
the difference $\bar S\setminus S$ by $\partial S$. The set $\pi(Q)$ is the preimage of $S$ under the map $F:\CC^m\to\C[A_0]\oplus\ldots\oplus\C[A_k]$, introduced before Definition \ref{defbif0}, and the preimage of the closure $\bar S$ has pure codimension $c$ by the following two reasons.

1) We have $\codim F^{-1}(\bar S)\geqslant c$, because $F^{-1}(\bar S)$ is the union of $\pi\{f_0^v=\ldots=f_k^v=0\}$ over all linear functions $v:\Z^n\to\Z$, and $\codim \pi\{f_0^v=\ldots=f_k^v=0\}$ is estimated by Proposition \ref{irrelpure}.1.

2) We have $\codim F^{-1}(\bar S)\leqslant c$, because $\codim \bar S=c$ by Lemma \ref{lunivres}.1.

Thus, the closure $\overline{\pi(Q)}$ is strictly smaller than the preimage of the closure $\bar S$, because a component $R'$ of the pure codimension $c$ set $F^{-1}(\bar S)$ contains the component $R\subset\overline{\pi(Q)}$ of higher codimension.

Consider a small perturbation $\tilde f_0,\ldots,\tilde f_k$ of $f_0,\ldots,f_k$, the corresponding perturbation $\tilde F$ of $F$, and denote $\{\tilde f_0=\ldots=\tilde f_k=0\}$ by $\tilde Q$. For generic $\tilde F$, its image intersects $\bar S$ and the closure of $\partial S$ properly (\cite{detandnewt}), thus the closure of $\pi(\tilde Q)$ equals $\tilde F^{-1}(\bar S)$. 
We summarize:

-- The intersection numbers of $\pi(Q)$ and $\pi(\tilde Q)$ with $L$ (counting multiplicities) are equal (by the result of Step III).

-- The intersection numbers of $\tilde F^{-1}(\bar S)$ and $F^{-1}(\bar S)$ with $L$ are equal
(apply (1) to $\pi\{f_0^v=\ldots=f_k^v=0\}$ for every $v$, and recall that $F^{-1}(\bar S)$ is their union).

-- The set $\tilde F^{-1}(\bar S)$ equals the closure of $\pi(\tilde Q)$, while $F^{-1}(\bar S)$ contains both $\pi(Q)$ and the codimension $c$
set $R'$, not contained in the closure of $\pi(Q)$.

Thus, denoting the intersection number (counting multiplicities) by ``$\circ$'', we have:
\begin{align*}
R'\circ L\leqslant F^{-1}(\bar S)\circ L - \pi(Q)\circ L= \\
=\tilde F^{-1}(\bar S)\circ L - \pi(Q)\circ L = \pi(\tilde Q)\circ L - \pi(Q)\circ L =0,
\end{align*}
i.e. the intersection number of $R'$ with a generic plane of complementary dimension equals 0, which is a contradiction. Thus, all components of the closure of $\pi(Q)$ have the same codimension $c$. $\quad\Box$

By Proposition \ref{irrelpure}.2, for $\min\dim A<-1$, the bifurcation set of $\pi$ equals the closure of $\pi(Q)$, which has codimension greater than 1, and whose tropicalization was described in \cite{st}, \cite{kaz}, etc. (it equals the projection of the tropical intersection of the dual complexes to the convex hulls of $H_0,\ldots,H_k$). Thus, the case $\min\dim A<-1$ is not interesting in the framework of our study of bifurcation sets, and we always assume that $A$ is relevant in what follows.
\par

\subsection{General case for relevant $A$}\label{Sgencaserel} We now choose a smooth toric compactification $T$ of the torus $\CC^n$, whose fan $\Gamma$ is compatible with the convex hulls of $A_0,\ldots,A_k$. Then every $f_i$ defines the section $\tilde f_i$ of the $A_i$-line bundle on $T\times\CC^m$ (see e.g. \cite{kh77} for
the construction of the line bundle on a toric variety, corresponding to a polytope). We denote the codimension 1 orbits of $T$ by $T_j,\, j=1,2,\ldots$, and consider the Bertini discriminant $\bar B_{\bar\pi}$ (Definition \ref{defbertdiscr}) for the projection $\bar\pi$ of $\bar Q=\{\tilde f_0=\ldots\tilde f_k=0\}$ to $\CC^m$, endowed with the divisors $D_j=(T_j\times\CC^m)\cap\bar Q$.
\begin{theor} \label{thmain} 1) If $\min\dim A\geqslant -1$ and $f$ is $0$-nondegenerate, then the Bertini discriminant $\bar B_{\bar\pi}$ is a hypersurface. \newline 2) If, moreover, $A$ is relevant and $f$ is $1$-nondegenerate, then we have $E_\pi=F\circ E_A$ and $|E_\pi|=B_\pi=\bar B_{\bar\pi}$ (in particular,
all of them are hypersurfaces), and the Newton polytope of $E_\pi$ equals $$\sum_{a_0+\ldots+a_k=n+1}\MP(\underbrace{H_0,\ldots,H_0}_{a_0},\ldots,\underbrace{H_k,\ldots,H_k}_{a_k}).$$
\end{theor}
See also \cite{gu} for the computation of the Euler characteristic of the bifurcation set $B_\pi$ in some simple cases.
\begin{lemma} \label{smoothnondeg} If $f$ is $c$-nondegenerate, then $\bar\pi$ is $c$-smooth (Definition \ref{defbertdiscr}).
\end{lemma}
\textsc{Proof.} We should prove that $\codim\bar\pi(\bigcup_I \sing Q_I)>c$. For every $I$, the set $Q_I$ has the form $\{f^u=0\}$ for some linear function $u:\Z^n\oplus\Z^m\to\Z$
that vanish on $\{0\}\times\Z^m$. Applying Lemma \ref{tropcover} with $\Sigma=\sing Q_I$, we have $\codim\pi\sing Q_I\geqslant c$, otherwise
the tropical fan of $Q_I$ contains a $c$-dimensional linear function $v$, which means that $\{f^v=0\}$ is not smooth,
i.e. $f$ is not $c$-nondegenerate. $\quad\Box$

\noindent\textsc{Proof of Theorem \ref{thmain}.} By Theorem \ref{thgeneral}.2 (for $\pi$ and $F$
as described here and $p$ as described in Subsection \ref{Sunivcase}), we have $\bar B_{\bar\pi}=F^{-1}(\bar B_{\bar p})$.
By Theorem \ref{lozung0}, $\bar B_{\bar p}$ is a hypersurface, thus $\bar B_{\bar\pi}$ is of codimension 1 or 0,
and it remains to exclude the possibility of zero codimension.
By definition, the set $\bar B_{\bar\pi}$ is the union of the projections of the sets \begin{center} $\Sigma_v=\sing\{f^v=0\}\cup\Bigl\{$
critical points of the projection of $\{f^v=0\}$ to $\CC^m\Bigr\}$ \end{center} over all linear functions $v:\Z^n\oplus\Z^m\to\Z$
that vanish on $\{0\}\times\Z^m$. Applying Lemma \ref{tropcover} to $\Sigma=\Sigma_v$,
we have $\codim\pi\Sigma_v>0$ for every $v$, otherwise the tropical fan of $\Sigma_v$ contains a 0-dimensional linear function
$u$, which means that $\Sigma_u$ is not empty, i.e. the projection of $\{f^u=0\}$ to $\CC^m$
is not a locally trivial fibration of smooth varieties. Since this projection is a locally trivial fibration
because of $\dim u=0$, we conclude that $\{f^u=0\}$ is not smooth for some 0-dimensional $u$, which contradicts 0-smoothness.

If, moreover, $\bar\pi$ is 1-smooth by Lemma \ref{smoothnondeg}, we can apply Theorem \ref{thgeneral}.2 and 3 with $\pi$ and $F$
as described here and $p$ as described in Subsection \ref{Sunivcase}. We get  
$E_\pi=F\circ E_A$ and $\bar B_{\bar \pi}=F^{-1}\bar B_{\bar p}$. Since the $\bar B_{\bar p}=|E_A|$ by Theorem \ref{thmainuniv},
the desired equality follow from Theorem \ref{thgeneral}.1.

\textsc{Newton polytope.} Choose finite subsets $B_1,\ldots,B_{m-1}\subset\Z^m$ and generic polynomials $g_i\in\C[B_i]$,
and compute the Euler characteristic of the set $f_0=\ldots=f_k=g_1=\ldots=g_{m-1}=0$ in the following two ways.

I. Let us count the desired Euler characteristic fiberwise, for every fiber of the projection $\{f_0=\ldots=f_k=g_1=\ldots=g_{m-1}=0\}\to\{g_1=\ldots=g_{m-1}=0\}$. Then, by the characterization of the Euler discriminant, obtained on Step II of the proof of Theorem \ref{thbasech}, the Euler characteristic equals $$e\{g_1=\ldots=g_{m-1}=0\}\cdot e\{ \mbox{ generic fiber of
} \pi\}+\{g_1=\ldots=g_{m-1}=0\}\circ E_\pi$$ (recall that $\circ$ is the intersection number counting multiplicities in $\CC^m$). Denoting the Newton polytope of $E_\pi$ by $X$, counting both Euler characteristics by the Khovanskii formula \cite{kh77} and the intersection number by the Kouchnirenko-Bernstein formula \cite{bernst}, we have (up to the sign $(-1)^{n-k}$) $$\MV(B_1,\ldots,B_{m-1},B_1+\ldots+B_{m-1})\cdot\sum_{a_0+\ldots+a_k=n}\MV(\underbrace{A_0,\ldots,A_0}_{a_0},\ldots,\underbrace{A_k,\ldots,A_k}_{a_k})+$$ $$+\MV(X,B_1,\ldots,B_{m-1}).$$
We allow a slight abuse of notation here, writing the mixed volume of finite sets instead of the one for their convex hulls.

II. Applying the Khovanskii formula \cite{kh77} to the system $f_0=\ldots=f_k=g_1=\ldots=g_{m-1}=0$, and removing the terms that are zero by Lemma \ref{lequiv1},
we have (up to the sign $(-1)^{n-k}$): $$\sum_{a_0+\ldots+a_k=n}\MV(B_1,\ldots,B_{m-1},B_1+\ldots+B_{m-1},\underbrace{H_0,\ldots,H_0}_{a_0},\ldots,\underbrace{H_k,\ldots,H_k}_{a_k})+$$ $$+\sum_{a_0+\ldots+a_k=n+1}\MV(B_1,\ldots,B_{m-1},\underbrace{H_0,\ldots,H_0}_{a_0},\ldots,\underbrace{H_k,\ldots,H_k}_{a_k}).$$
The first term equals $$\MV(B_1,\ldots,B_{m-1},B_1+\ldots+B_{m-1})\cdot\sum_{a_0+\ldots+a_k=n}\MV(\underbrace{A_0,\ldots,A_0}_{a_0},\ldots,\underbrace{A_k,\ldots,A_k}_{a_k})$$ by Lemma \ref{mvprod}.
Comparing (I) and (II), we have:
$$\MV(X,B_1,\ldots,B_{m-1})=\sum_{a_0+\ldots+a_k=n+1}\MV(B_1,\ldots,B_{m-1},\underbrace{H_0,\ldots,H_0}_{a_0},\ldots,\underbrace{H_k,\ldots,H_k}_{a_k}).$$
By Lemma \ref{lekh}, the desired sum of mixed fiber polytopes satisfies this equation for all $B_1,\ldots,B_{m-1}$,
and is uniquely determined by it.$\quad\Box$.

\begin{utver}\label{nonempty1} Assume that $\min\dim A=-c$ and $f$ is $c$-nondegenerate. \newline
1) Assume, moreover, that $c>0$. 
The bifurcation set $B_\pi$ is empty if, for some non-empty essential subtuple $A_{i_1},\ldots,A_{i_p}$ of $A$ (see Definition \ref{defess}), the projection
$H_{i_1}+\ldots+H_{i_p}\to\R^n$ is injective, and is non-empty of pure codimension $c$ otherwise. \newline
2) Assume that $c=0$ (i.e. $A$ is linearly independent). The bifurcation set is empty if the projection
$H_0+\ldots+H_k\to\R^n$ is injective, and is a non-empty hypersurface otherwise.
\end{utver}
Note that the condition for emptiness of the bifurcation set of $\{f=0\}\to\CC^m$ in the case $c>0$ implies that $\{f=0\}$ itself is empty.

\textsc{Proof.} The codimension is computed in Proposition \ref{irrelpure}.2 and Theorem \ref{thmain}.2 for $c>0$ and $c=0$ respectively. To check emptiness for $c\leqslant 1$, apply the criterion of triviality of the mixed fiber polytope (Theorem \ref{trivmp}) to the Newton polytope of the discriminant, computed in Theorems \ref{thmain}.2. To reduce the case $c>1$
to $c=1$, take the standard simplex $C\subset\Z^{c-1}$, and consider $A'_i= C\times A_i$, $H'_i=C\times H_i\in\Z^{c-1}\oplus\Z^n\oplus\Z^m$ and generic $f'_i\in\C[H'_i]$. Then $\min\dim A'=-1$, and the desired emptiness assertion for $f$ follows from the same for $f'$. $\quad\Box$

\end{document}